\documentclass{gSTA2e}

\theoremstyle{plain}
\newtheorem{thm}{Theorem}[section]

\theoremstyle{remark}
\newtheorem{remark}[thm]{Remark}
 
 \usepackage{amsmath}
\usepackage{amssymb}
\usepackage{epsfig}
\usepackage{graphicx}
\usepackage{upgreek}
 
\def\N{{\rm I\kern-.15em N}}
\def\R{{\rm I\kern-.2em R}}
\def\Z{{\rm Z\kern-.26em Z}}
\makeatother

\newcommand{\be}{\begin{eqnarray}}
\newcommand{\ee}{\end{eqnarray}}
\newcommand{\bq}{\begin{eqnarray*}}
\newcommand{\eq}{\end{eqnarray*}}

\newcommand{\eps}{\varepsilon}

\begin{document}

\title{Tests for Time Series of Counts Based on the Probability Generating Function}
 
\author{\v{S}\'arka Hudecov\'a$^{a}$  and  Marie Hu\v{s}kov\'a$^{a}$    and   Simos G. Meintanis$^{bc}$$^{\ast}$\thanks{$^\ast$Corresponding author. Email: simosmei@econ.uoa.gr
\vspace{6pt}}\\\vspace{6pt}  $^{a}${\em{Charles University in Prague, Faculty of Mathematics and Physics, Department of
Probability and Mathematical Statistics,
            Sokolovsk\'a 83, CZ\,--\,186\,75 Praha 8, Czech Republic}};\\
$^{b}${\em{National and Kapodistrian University of Athens,
Department of Economics, 1 Sophocleous Str., and Aristidou Str., 105 59 Athens,
Greece}}\\
$^{c}${\em{Unit for Business Mathematics and Informatics, North West University, Potchefstroom,
South Africa$^{\dagger}$\thanks{$^{\dagger}$On sabbatical leave from the University of Athens}}}
\\\vspace{6pt}\received{received November 2013} }

\maketitle

\begin{abstract}
We propose testing procedures for the hypothesis that a given set
 of discrete observations may be formulated as a particular time series of counts with a specific
 conditional law. The new test statistics incorporate  the empirical probability generating function
 computed from the observations. Special emphasis is given to the popular models of integer
 autoregression and
 Poisson autoregression. The asymptotic properties of the proposed
 test statistics are studied under the null hypothesis as well as under alternatives. A Monte Carlo
  power study on bootstrap versions of the new methods is included as well as real--data examples.


\begin{classcode} 62M10, 62F03, 62F40  \end{classcode}

\end{abstract}

\section{Introduction}\label{sec_0}

Let $\{Y_t\}_{t\in \mathbb N}$ denote a time series of counts for which, conditionally on the past, the corresponding stochastic structure is fully specified by a family of laws indexed by a certain parameter. Such models include the model of integer autoregression (INAR) as well as the integer autoregressive conditionally heteroscedastic (INARCH) model. These two models have received enormous attention in the past as they are known to fit empirical data in diverse areas of application. The objective here is to construct goodness--of--fit (GOF) statistics for distributional assumptions regarding these count time series models.
In the classical (continuous type) framework of time--series models, this aspect of modelling has drawn considerable attention recently; see
\cite{koul2013, koul-mimoto, koul-surgailis, horvath2006,horvath2004}.
The standard approach in constructing GOF tests is to estimate the corresponding density  or distribution function and thereby construct versions of the Kolmogorov--Smirnov, Cram\'er--von Mises and Bickel--Rosenblatt statistics. However, there is also the alternative route of employing empirical transforms, such as the empirical Laplace transform and the empirical characteristic function, to the same problem. This idea was first put forward by Epps \cite{epps} and has recently been followed by
Cuesta Albertos et al. \cite{cuesta}, Ghosh \cite{ghosh}, and Klar et al. \cite{klar}.
 
Now turning to count data, Fokianos and Neumann \cite{fokianos-neumann} have considered GOF tests for the regression function in parametric versions of count time series. Here we consider another aspect of such series, namely the aspect of correctly specifying the conditional distribution of observations. In doing so we employ the aforementioned approach of empirical transforms and use marginal quantities by integrating up with respect to the conditioning argument in a spirit analogous to that in \cite{rothe}. Specifically, the test statistics comes in the form of a weighted L2--type distance between a nonparametric estimate of the marginal probability generating function of the observations and a semiparametric estimate of the same quantity imposed by the model.   Recall that if $Y$ is an arbitrary integer--valued random variable then its   proba\-bi\-lity generating function (PGF) is defined as
 $g_Y(u)=\mathbb E[u^Y], \ u \in [0,1]$.

The remainder of this paper runs as follows. In Section 2 we specify the time series and the GOF problems considered, and construct the corresponding test statistics for the first order models. Computational issues are addressed in Section~3. In Section 4 the asymptotic properties of these statistics are studied both under the null hypothesis as well as under alternatives. Corresponding proofs are postponed to the appendix. In Section 5 we propose estimators for the parameters involved in our test statistics and suggest bootstrap versions of the new tests. Possible generalization of our approach to models of higher order is described in Section 6. In Section 7 we report the results on a Monte Carlo study, and the article concludes with some real data examples in Section 8.

\section{Test statistics}\label{sec_1}

Let ${\cal{I}}_{t}$ be  the information set available at time $t$, i.e., ${\cal{I}}_{t}=\sigma\{Y_s,\ s\leq t\}$ is a $\sigma$-field generated by the past values of the  series $\{Y_t\}_{t\in \mathbb N}$.   We assume that the conditional distribution of $Y_t$ given ${\cal{I}}_{t-1}$ can be described using a specific (cumulative) distribution function $F$
depending  on ${\cal{I}}_{t-1}$ and also on an unknown parameter $\vartheta$
 in the following way:
\begin{equation} \label{basicINAR}
Y_t | \: {\cal{I}}_{t-1} \sim F(\cdot; {\cal{I}}_{t-1},\vartheta),
\end{equation}
where $X \sim F(\cdot;\cdot,\cdot)$ is interpreted as `the random variable $X$ has $F(\cdot;\cdot,\cdot)$ as its distribution function'.
 
There exist two main  classes of models that may be formulated in this way: The INAR model and INARCH model, both admitting several specifications and generalizations. We start with the most basic formulations.

\noindent {\bf{(i)}} INAR model: For the INAR(1) model (see \cite{mckenzie1985,al-osh-1987,alzaid1988})
eqn. \eqref{basicINAR} is specified as
\begin{equation} \label{modelF(2.2}
 Y_t | \: {\cal{I}}_{t-1} \sim {\cal{B}}(Y_{t-1},p) \ast G_{\varepsilon},
 \end{equation}
 where   $\ast$ denotes convolution, ${\cal{B}}(\nu,p) $ is a
 binomial distribution with parameter $(\nu,p)$  and
 $G_{\varepsilon}$  is a distribution.{More specifically the INAR(1) model can be
constructed as
\begin{equation}\label{modelINAR} Y_t=\sum_{i=1}^{Y_{t-1}}U_{t,i}
+\varepsilon_t,\quad t=1,\ldots,
\end{equation} where $U_{1,1},\ldots$ are independent binary random variables such that
$P(U_{t,i}=1)=p,\, P(U_{t,i}=0)=1-p,$ $p\in(0,1)$, and such that $U_{t,1}, U_{t,2},\dots$ are independent of $Y_{t-1}$ and $\varepsilon_t$. The  variables
$\varepsilon_t,\, t=1,\ldots $ are i.i.d. discrete nonnegative with distribution function $G_{\varepsilon}$ and PGF $g_{\varepsilon}(u)$ with finite  variance  
such that $\varepsilon_t$ is independent of
  $Y_{t-1}$.      
 
Generalizations of the INAR(1) to higher order models were proposed, see \cite{alzaid1990,du-li},
    while a related review article is \cite{mckenzie2003}.
 
 Specific instances of INAR(1) result when $G_{\varepsilon}$ is known to belong to a family of
  distributions ${\cal{G}}_\Theta$, the most popular case being when ${\cal{G}}_\Theta$ is
  the Poisson family, in which case $\vartheta=(p,\theta)'$ with $\theta \in \Theta=(0,\infty)$ denoting the Poisson parameter. In fact, the pro\-perty that  the marginal distribution of the observations is from the same family as the distribution of innovations $\varepsilon_t$ characterizes the Poisson law in the context of INAR mo\-dels. Wei$\upbeta$ \cite{weiss} however points out that both the Negative Binomial as well as Consul's generalized Poisson distribution may well serve as marginal laws under INAR mo\-dels, particularly in view of overdispersed data. (In fact the set of all possible marginal laws coincides with the family of discrete self--decomposable distributions). Under the same type of data Pavlopoulos and Karlis \cite{pavlopoulos-karlis} suggest a Poisson mixture for the innovation distribution, while Barczy et al. \cite{barczy} entertain the idea of an   INAR model with Poisson innovations, but containing outliers in the innovation distribution. Hence, and since under a fixed INAR model the law of the innovations uniquely determines the conditional as well as the marginal law of the observations, there is a clear need for GOF procedures in order to correctly identify the innovation distribution $G_\varepsilon$.

\medskip

\noindent {\bf{(ii)}} INARCH model: The INARCH(1) model is specified
by (\ref{basicINAR}) with
\begin{equation} \label{inarch}
 Y_t | \: {\cal{I}}_{t-1} \sim  \mathsf{Po}( \lambda_t),\quad \lambda_t=r(Y_{t-1}; \vartheta)
 \end{equation}
where  $\mathsf{Po}( \lambda_t)$ denotes Poisson distribution with parameter $\lambda_t$ , and $r(\cdot;\cdot)$ belongs to
 a specific parametric family of functions $\{r(Y_{t-1}; \vartheta),\,  \vartheta\in \Theta\})= \mathcal{G}_{\Theta}$.

A specific instance
is the Poisson linear autoregression; see \cite{fokianos2009}. Generalizations of the INARCH(1) model may be found in \cite{fok-tjos}, while \cite{fokianos2011} contains a nice synopsis of
this model as well as related models. Although here too the Poisson assumption is by far the most popular specification, alternative distributional specifications in eqn.~\eqref{inarch} such as the Negative Binomial INARCH model of Zhu \cite{zhu2011} and the INARCH with interventions of Fokianos and Fried \cite{fokianos-fried},  have also been proposed.
 
\subsection{{\rm{{\bf{INAR model}}}}}

We begin our GOF discussion with the INAR(1) model specified by
eqn. (\ref{modelINAR}).  Given the data $Y_t, \
t=1,...,T$, one wishes to test the null hypothesis
\medskip

\noindent $ {\cal{H}}_{01}$:   $\{Y_t\}_{t\in \mathbb N}$ follow model (\ref{modelINAR})
for some $p\in (0,1)$ and   some PGF $g_{\varepsilon} $ belonging to a
given family ${\cal{G}}_\Theta=\{ g_{\varepsilon}(\cdot;\theta);\, \theta\in
\Theta\}$  with $\Theta$ being an open subset of $\mathbb{R}$,

\medskip
\noindent against a general alternative that  $ {\cal{H}}_{01}$ does not hold. Notice that  by results in  \cite{du-li}   the sequence $\{Y_t\}_{t\in \mathbb N}$ is stationary  and ergodic under $ {\cal{H}}_{01}$. Specifically in \cite{du-li} it is shown that for $p<1$, there exist a unique stationary and ergodic solution satisfying eqn. \eqref{modelINAR}, which is produced by the Markov chain in this equation.  \cite{du-li} also provide conditions for a stationary  and ergodic INAR model of arbitrary order.

 We suggest to test the null hypothesis ${\cal{H}}_{01}$ by
means of the test statistic
\begin{equation}\label{testINAR}
 S_T=T
\int_0^1 \left(\widehat g_{T,Y}(u)-\widehat g_{T0}(u)\right)^2
w(u)du,
\end{equation}
where $w(u)$ is a nonnegative weight function,
\begin{equation} \label{eq26} \widehat
g_{T,Y}(u)=\frac{1}{T} \sum_{t=1}^{T} u^{Y_t},
\end{equation}
is the
empirical PGF which is a non--parametric estimate of the PGF of
$\{Y_t\}$, while  $ \widehat g_{T0}(u)$ will be a  {semiparametric}
estimate of the same PGF under the model specified by the null
hypothesis ${\cal{H}}_{01}$.

 Notice that  in model (\ref{modelINAR})   the PGF of
$\sum_{i=1}^{N}U_{t,i}$ is given by $(1+p(u-1))^N$ and that for the marginal population PGF of $Y_t$  $g_{Y_t}(\cdot)$  we have
\begin{align}\label{eq27a}
 g_{Y_t}(u)=&\mathbb E\left[u^{Y_t}\right]=\mathbb E \left [\mathbb E(u^{Y_t}| Y_{t-1})\right]\\
= \nonumber &\mathbb E \left  [ (1+p(u-1))^{Y_{t-1}}g_{\varepsilon_t} (u
) \right]
 =g_{\varepsilon_t}(u) g_{Y_{t-1}}( 1+p(u-1)), \quad u\in [0,1],
\end{align}
  for some $p\in(0,1)$ and some $\theta\in \Theta$.  Since under the considered assumptions $\{Y_t\}$ are  strictly stationary and $\{\varepsilon_t\}$ are i.i.d. random variables  we can  write $g_{Y}(\cdot)$ and $g_{\varepsilon}(\cdot)$  instead of
$g_{Y_t}(\cdot)$ and $g_{\varepsilon_t}(\cdot)$, respectively.
 Under the null
hypothesis ${\cal{H}}_{01}$  the relation (\ref{eq27a})  reduces to
\begin{equation}\label{eq27}
 g_{Y}(u)=
 g_{
 \varepsilon}(u; \theta) g_{Y}( 1+p(u-1)), \quad u\in [0,1],
\end{equation}
where  $g_{\varepsilon}(\cdot; \theta)$ denotes the PGF of  $\{\varepsilon_t\}$ under the null hypothesis.  Recently \cite{MK2014} employ a non--parametric against a fully parametric estimate of the PGF of the observations in testing within the INAR context. Here however we follow a different approach.} Assume that $\widehat \vartheta_T=(\widehat p_T,\widehat\theta_T)'$ is a suitable estimator of $\vartheta$ constructed from $Y_1,\dots,Y_T$.
Then based on eqn.~\eqref{eq27}, 
  a natural semiparametric estimate of the marginal PGF is
\begin{equation} \label{estimnull}
\widehat g_{T0}(u)=g_{\varepsilon}(u;\widehat \theta_T)   \widehat
g_{T,Y}(1+\widehat p_T(u-1)), \quad u\in [0,1],
\end{equation}
where $g_{\varepsilon}(u;\widehat \theta_T)$  and $\widehat{g}_{T, Y}(1+\widehat p_T(u-1))$
are the PGFs  of
$\varepsilon_t$  under the null hypothesis with $\theta$ replaced by an estimator  $\widehat
\theta_T$  and  of $\widehat g_{T,Y}(u)$ defined in
eqn.~\eqref{eq26} computed at the point $1+\widehat p_T(u-1)$.
 
\subsection{{\rm{{\bf{INARCH model}}}}}

Likewise, for the INARCH model one wishes to test an analogous GOF null hypothesis which may be formulated as

\medskip

\noindent $ {\cal{H}}_{02}$:   $\{Y_t\}_{t\in \mathbb N}$ follow model
$$ 
 Y_t | \: {\cal{I}}_{t-1} \sim  \mathsf{Po}( \lambda_t),\quad \lambda_t= \theta_1+Y_{t-1} \theta_2
 $$
  with $(\theta_1, \theta_2)'\in \Theta= \{ (\theta_1, \theta_2)'; \theta_1>0, \theta_2\in [0,1)\}$.
 
 \medskip

 Notice that under $ {\cal{H}}_{02}$ the sequence   $\{Y_t\}_{t\in \mathbb N}$ is stationary and ergodic. In fact, \cite{ferland} prove strict stationarity of the more general INGARCH model of arbitrary orders under assumptions entirely parallel to those of an ARMA model. Also \cite{fokianos2009} use a perturbed INGARCH model which can be made arbitrarily close to the corresponding (unperturbed) INGARCH model, in order to obtain ergodicity properties and to prove the asymptotic properties of estimators of the parameters in the latter model. Stationarity and ergodicity properties of the Poisson INGARCH model with non--linearly specified intensity $\lambda_t$ are discussed by \cite{Neumann2011}.

For the null hypothesis  ${\cal{H}}_{02}$ we employ again
the test statistic in eqn.~\eqref{testINAR}, but one needs to find a semiparametric
estimate of the PGF reflecting now the null hypothesis
${\cal{H}}_{02}$. To this end, we first compute the corresponding
marginal population PGF as follows:
 \begin{equation} \label{eq29} g_{Y_t}(u)=\mathbb
E\left[u^{Y_t}\right]=\mathbb E \left [\mathbb E(u^{Y_t}|
Y_{t-1})\right]=\mathbb E \left  [g_{Y_t}(u;\theta_1+\theta_2
Y_{t-1})\right], \quad u\in [0,1].
\end{equation}

Then a natural semiparametric estimate of the marginal PGF should be based on an estimate of $\mathbb E \left  [g_{Y_t}(u;\theta_1+\theta_2 Y_{t-1})\right]$, where $(\theta_1,\theta_2)'$ is  replaced by a suitable estimate $(\widehat \theta_{1T},\widehat \theta_{2T})'$.

In case of Poisson conditionals as in eqn.~\eqref{inarch}, $\vartheta:=(\theta_1,\theta_2)\in \Theta \subset \mathbb{R}^2$ with $\Theta=(0,\infty)\times(0,1)$. Recall also that the PGF of the Poisson
distribution with mean $\theta$ is given by $e^{\theta(u-1)}$.  Then we have from eqn.~\eqref{eq29}
\begin{align}\label{INARCHa}
g_Y(u)=&\mathbb E
\left  [e^{(\theta_1+\theta_2 Y_{t-1})(u-1)}\right]\\
\nonumber=&e^{\theta_1(u-1)}\mathbb E
\left  [e^{\theta_2(u-1)Y_{t-1}}\right]=e^{\theta_1(u-1)} g_{Y}\left(e^{\theta_2(u-1)}\right).
\end{align}
 In view of eqn.~\eqref{INARCHa}, a semiparametric estimate of the marginal PGF under $\mathcal{H}_{02}$ is
\be
 \widehat g_{T0}(u)=e^{\widehat\theta_{1T}(u-1)} \widehat g_{T,Y}(e^{\widehat\theta_{2T}(u-1)}),
\ee
where $\widehat{g}_{T,Y}(e^{\widehat\theta_{2T}(u-1)})$ is the empirical PGF in eqn.~\eqref{eq26} computed at the point $e^{\widehat \theta_{2T}(u-1)}$.

\section{Computations}\label{sec2}

In this subsection we simply write $\widehat \vartheta$ instead of $\widehat \vartheta_T$ for the estimator. For the INAR(1) model in eqn.~\eqref{modelINAR} we have from eqns.~\eqref{testINAR} and \eqref{estimnull} by straightforward algebra
\be
S_T=\frac{1}{T}\sum_{t,s=1}^T \left\{I^{(0)}(Y_{t,s},0)+I^{(2)}(0,Y_{t,s})-2I^{(1)}(Y_t,Y_s)\right\},
\ee
where $Y_{t,s}=Y_t+Y_s$, and
\begin{equation}\label{eq32}
I^{(m)}(x,y):=I_{w,g}^{(m)}(x,y)=\int_0^1 (g_\varepsilon(u;\widehat{\theta}))^m u^x \: (1+\widehat{p}(u-1))^y \: w(u)du.
\end{equation}
\noindent To proceed any further we will need to assume a specific family ${\cal{G}}_\Theta$ under the null hypothesis ${\cal {H}}_{01}$ and fix the weight function $w(\cdot)$. In particular if we let ${\cal{G}}_\Theta$ be the Poisson family of distributions (so that $g_\varepsilon(u;\theta)=e^{\theta(u-1)}$), and choose $w(u)=u^a, \ a\geq 0$, as weight function, we have from eqn.~\eqref{eq32}
\[
I^{(0)}(x,0)=\int_0^1 u^x u^a du=\frac{1}{1+a+x},
\]
\[
I^{(2)}(0,x)=\int_0^1 e^{2\widehat{\theta}(u-1)} U_{\widehat{p}}^x u^a du=e^{-2\widehat{\theta}}\sum_{\ell=0}^x \sum_{m=0}^\ell (-1)^{\ell-m}\left(\begin{array}{c}x\\\ell\end{array}\right) \left(\begin{array}{c}\ell \\ m\end{array}\right) \widehat{p}^\ell  J(a+m,2\widehat{\theta}),
\]
and
\[
I^{(1)}(x,y)=\int_0^1 e^{\widehat{\theta}(u-1)} u^x U_{\widehat{p}}^y u^a du=e^{-\widehat{\theta}}\sum_{\ell=0}^y \sum_{m=0}^\ell (-1)^{\ell-m} \left(\begin{array}{c}y\\\ell\end{array}\right) \left(\begin{array}{c}\ell \\ m\end{array}\right) \widehat{p}^\ell J(a+x+m,\widehat{\theta}),
\]
where we have used the notation $U_{\widehat{p}}=1+\widehat{p}(u-1)$ and
\[
J(\lambda,\mu)=\int_0^1 u^\lambda e^{\mu u} du=\sum_{k=0}^\infty \frac{1}{1+\lambda+k} \:\frac{\mu^k}{k!}=\frac{(-\mu)^{-\lambda}}{\mu} \left[\Gamma(1+\lambda,-\mu)-\Gamma(1+\lambda)\right].
\]

\smallskip

Likewise for the INARCH(1) model in eqn.~\eqref{inarch} with a Poisson conditional distribution we have from eqns.~\eqref{testINAR} and \eqref{INARCHa} by straightforward algebra
\be S_T=\frac{1}{T}\sum_{t,s=1}^T \left\{ I^{(0)}(Y_{t,s},0)+e^{-w_{ts}}J(a,w_{ts})-2e^{-w_{s}}J\left(a+Y_t,w_s\right)\right\},
\ee
where $w_{ts}=2\widehat{\theta}_1+\widehat{\theta}_2 Y_{t,s}$ and $w_{s}=\widehat{\theta}_1+\widehat{\theta}_2 Y_{s}$.
 
 \section{Asymptotic results}\label{sec3}

  Here we study the limit behavior ($T\to \infty$) of the test statistic $S_T$ both for INAR(1) and INARCH(1) under  the null hypothesis as well as  under alternatives.  In what follows, we present the results for the INAR(1) in detail. Corresponding results for the  INARCH(1) are derived in an analogous manner and therefore are presented with less detail.
 
\subsection{{\rm{{\bf{INAR(1) model}}}}}

    Recall that for the INAR(1)  model   formulated in (\ref{modelINAR})    the PGFs of $Y_t$ and $\varepsilon_t$ satisfy (\ref{eq27a})
and  under
${\cal{H}}_{01}$
\[
  g_{\varepsilon}(u)\in\mathcal{G}_{\Theta}=\{g_{\varepsilon}(u; \theta), \, \theta\in
  \Theta\},\quad u\in  [0,1],
\]
 i.e.  $g_{\varepsilon}(\cdot; \theta)$ is specified up to a parameter
 $\theta\in \Theta$, with $\Theta$ being an open set in $\mathbb{R}$  and (\ref{eq27})  holds true.

 \smallskip

Denote the true value of $\vartheta=(p,
 \theta)' $ under the null hypothesis ${\cal{H}}_{01}$  by $\vartheta_0=(p_0,\theta_0)'$.
 To study the limit distribution under the null hypothesis ${\cal{H}}_{01}$  we assume the following:

 \begin{enumerate}

 \item[(A.1)] $\{Y_t\}_{t\in \mathbb N}$ is a sequence of random variables
 (\ref{modelINAR}) with  $\{\varepsilon_t\}_{t\in \mathbb N}$  being  a sequence of i.i.d. discrete nonnegative random variables with  finite second moment and PGF
 $g_{\varepsilon} (\cdot;\theta),\, \theta\in \Theta$, where $\Theta$ is an
 open subset of  $\mathbb{R}$.

 \item[(A.2)] $g_{\varepsilon}(u; \theta)$  has  the  first partial derivative w.r.t $\theta$
 for all $u\in[0,1]$  fulfilling  Lipschitz condition:
     \[
    \left| \frac{\partial g_{\varepsilon}(u;\theta)}{\partial \theta}-
    \frac{\partial g_{\varepsilon}(u;\theta)}{\partial \theta}\big{|}_{\theta=\theta_0}\right|\leq D_1 |\theta_0
    -\theta| v (u),\quad u\in[0,1],\quad     |\theta-\theta_0|\leq D_2
    \]
    and
     \[
    \left|\frac{\partial g_{\varepsilon}(u;\theta)}{\partial \theta}\right|\leq D_3  v(u),\quad
    u\in[0,1],\quad   |\theta-\theta_0|\leq D_2
    \]
     for some $D_j>0,\, j=1,2,3$, and some  measurable function $v(\cdot)$.
 \item[(A.3)] $\quad
 0<\int_0^1 w(u) du <\infty, \quad  \int_0^1 w(u)v^2(u) du<\infty.$

 \item[(A.4)] $\widehat{\vartheta}_T=(\widehat p_T,\widehat \theta_T)'$ is estimator of the true value of  $\vartheta_0=(p_0, \theta_0)'$
  satisfying
  \[
  \sqrt{T} (\widehat \vartheta_T - \vartheta_0)=\frac{1}{\sqrt{T}} \sum_{t=1+q}^T \ell ({\bf{Y}}_{t-q}; \vartheta_0)+  o_P(1),
  \]
  where ${\bf{Y}}_{t-q}=(Y_t,...,Y_{t-q})'$, and $\ell=(\ell_1,\ell_2)'$ is such that
   $ \ell_{j} ({\bf{Y}}_{t-q};  \vartheta_0),\, j=1,2$ (with fixed $q \geq 1$), are assumed to be martingale difference sequences with    finite variances.
   \end{enumerate}

Define for $t\geq 2$
\begin{equation}\label{Zt}
Z_t(u;p,\theta)=u^{Y_t}-(1+p(u-1))^{Y_{t-1}} g_{\varepsilon} (u;\theta),
\end{equation}
and
\begin{equation}\label{calZt}
{\cal{Z}}_{t}(u;p, \theta)=Z_{t}(u;p,\theta)+h_1(u;p,\theta)\ell_1({\bf{Y}}_{t-q};p,\theta)
+h_2(u;p,\theta)\ell_2({\bf{Y}}_{t-q};p,\theta)
\end{equation}
for $ t\geq q+1$, where
  \begin{align*}\label{h1h2}
      h_1(u; p,\theta)=& \frac{\partial g_{Y}(1+p(u-1))}{\partial p} g_{\varepsilon} (u; \theta ),\notag\\
      h_2(u;p, \theta)=&  g_Y(1+p(u-1)) \frac{\partial g_{\varepsilon}(u;\theta)}{\partial \theta}.
      \end{align*}

In the following theorem we prove the main assertion on limit behavior of our test statistic $S_T$ under ${\cal{H}}_{01}$.

  \begin{thm}
  \label{thm_2_2} Let assumptions (A.1)-(A.4) be satisfied  in the model (\ref{modelINAR}). Then
 under the null hypothesis ${\cal{H}}_{01}$ the limit distribution ($T\to \infty$) of $S_T$ defined in eqn. (\ref{testINAR}) is the same as that of
\[
\int_0^1  \bar{\cal{Z}}_T^2(u; p_0, \theta_0) w(u) du,
\]
where
$
\bar{\cal {Z}}_T(u;p_0,\theta_0)=\frac{1}{\sqrt T} \sum_{t=q+1}^T{ \cal{Z}}_t (u;p_0,\theta_0)
$
with $\{ {\cal{Z}}_{t}(u;p, \theta) , \, u\in [0,1]\}$ as in eqn.~\eqref{calZt}. 
Moreover, the process $\{\bar{\cal{Z}}_T(u;p_0, \theta_0); \ u\in[0,1]\}$  converges in $ \mathcal{C}[0,1]$ to a zero--mean Gaussian process $\{{\cal{Z}}(u; p_0, \theta_0); \ u\in[0,1]\}$  with covariance structure $\mathbb {E}\Big[{\cal{Z}}_{q+1}(u_1;p_0,\theta_0){\cal{Z}}_{q+1}(u_2;p_0,\theta_0)\Big]$, $u_1, u_2 \in [0,1]$, 
and $S_T$ converges in distribution to
\[
\int_0^1  {\cal{Z}}^2(u; p_0, \theta_0) w(u) du.
\]
      \end{thm}

      \noindent The proof is postponed to the Appendix.

      \smallskip

\begin{remark}
 Note that there is no explicit form for the limit distribution function of the test statistic $S_T$, and that this distribution function depends on unknown quantities. Therefore, Theorem~\ref{thm_2_2} is not directly applicable for the purpose of approximating critical values and actually performing the test.  Nevertheless, when a consistent estimator of the covariance structure is available
  we can plug it instead of unknown quantities
   and the assertion of our theorem remains true.
     Alternatively, in Section 5, a properly chosen bootstrap  is proposed which provides an effective way for approximating the limit null distribution of $S_T$. In both cases of approximation however we need an estimator $\widehat{\vartheta}_T=(\widehat p_T,\widehat \theta_T)'$ of the parameters. In Section~\ref{sec5} we also construct these estimators.
\end{remark}

\bigskip

We now consider the behavior of the test statistic under alternatives of the type $g_{\varepsilon} \notin {\cal{G}}_\Theta$, which means that we still have model \eqref{modelINAR}  but the innovation distribution need not belong to the family ${\cal{G}}_\Theta$.

     We will assume that $(\widehat p_T,\widehat \theta_T)'$ has the property
     \begin{equation} \label{alternativeest}
      (\widehat p_T, \widehat \theta_T)' \stackrel{{\rm{P}}}{\rightarrow}  ( p,  \theta_A)'.
     \end{equation}
    for   $ p\in(0,1)$   being the true parameter value and  for $\theta_A  \in \Theta$.

     \begin{thm}
  \label{thm_2_3} 
     Let $\{Y_t\}_{t\in \mathbb N}$   follow the model (\ref{modelINAR}). Let   (\ref{alternativeest}) and (A.3) be  satisfied and let   also  $g_{\varepsilon}(u; \theta)$ be continuous in $\theta $ for all $u\in[0,1]$.    Then, as $T\to \infty$,
     \begin {equation} \label{lim-alternx}
     \frac{S_T}{T}  \stackrel{{\rm{P}}}{\rightarrow}  \int_0^1 \Big[ g_Y(1+ p(u-1))( g_{\varepsilon}(u)-g_{\varepsilon}(u; \theta_A))\Big]^2 w(u) du,
     \end{equation}
     where $p$ is the true value of the parameter.
     \end{thm}

 The proof is omitted since  it suffices to follow the line of the proof of Theorem 4.1 and  use stationarity and ergodicity of $\{Y_t\}_{t\in \mathbb N}$.

 \medskip

 However, the right--hand side of (\ref{lim-alternx}) is strictly positive unless the true innovation PGF
$g_{\varepsilon}(\cdot)$ coincides with the PGF $g_{\varepsilon}(\cdot; \theta_A)$ postulated under the null hypothesis ${\cal{H}}_{01}$.
 This and the uniqueness of the PGF implies  the consistency of the test which rejects the null hypothesis ${\cal{H}}_{01}$ for large values of the test statistic $S_T$ under fixed alternatives. We should point--out however, that despite the fact that the formulation of the alternative adopted here focuses exclusively on the innovation PGF, eqn. (7) clearly reflects not only this PGF but the entire INAR model as this model is specified by eqn.~\eqref{modelINAR}. Therefore our test is expected to also have non--negligible power against model violations. This feature of the test will be further illustrated by simulations in Section 7.

     It can be further proved that the
       test is also sensitive w.r.t.  local alternatives,   e.g., it is true if the innovation PGF may be written as $\widetilde g_{\varepsilon}(u)=g_{\varepsilon}(u)(1+(\kappa/\sqrt{T})f(u)), \ \kappa \neq 0$, where the function $f(\cdot)$ is such that $\widetilde g_{\varepsilon}(\cdot)$ is a PGF. The  derivation of the corresponding results however is quite technical and therefore we will not pursue this issue here any further.

       \subsection{{\rm{{\bf{INARCH(1) model}}}}}

 Now we turn to the INARCH(1) setup.
As already mentioned, the limit behavior of the test statistic $S_T$ for this setup can be obtained in a manner quite analogous  to the INAR(1) case. 

Denote the true values of $(\theta_1,
  \theta_2)'$ by $(\theta_{10},
  \theta_{20})'$ and assume also that $(\widehat \theta_{1T}, \widehat \theta_{2T})'$ are estimators of  $(\theta_1,
  \theta_2)'$ satisfying
 \begin{equation} \label{INARCHest}
        \sqrt T (\widehat \theta_{jT}-\theta_{j0})=\frac{1}{\sqrt T}\sum_{t=q+1}^T \ell_{j} ({\bf{Y}}_{t-q}; \theta_{10}, \theta_{20})+
        o_P(1), \quad j=1,2,
        \end{equation}
  where
   $ \ell_{j} ({\bf{Y}}_{t-q}; \theta_{10}, \theta_{20}),\, j=1,2$,
   for fixed  $q\geq 1$, are assumed to be martingale difference sequences with finite variances.

    \smallskip
     
Define for $t\geq 2$
\begin{equation}\label{eq45}
V_t(u;\theta_1, \theta_2)=u^{Y_t} -\exp\{(\theta_1+\theta_2Y_{t-1})(u-1)\},
\end{equation}
and 
\begin{equation}\label{calVt}
{\cal {V}}_{t}(u;\theta_{1},\theta_{2})=  V_{t}(u;\theta_{1},\theta_{2})+r_1(u;\theta_{1},\theta_{2})
\ell_1({\bf{Y}}_{t-q};\theta_{1},\theta_{2})
\\+r_2(u;\theta_{1},\theta_{2})
\ell_2({\bf{Y}}_{t-q};\theta_{1},\theta_{2})
\end{equation}
for $t\geq q+1$,
where
  \begin{align*}\label{r1r2}
      r_1(u;\theta_{1},\theta_{2})=& {\mathbb {E}} \Big[\exp\{(\theta_1+\theta_2Y_1)(u-1)\}(u-1)\Big],\notag\\
      r_2(u;\theta_1,\theta_2)=&  {\mathbb {E}} \Big[Y_{1} \exp\{(\theta_1+\theta_2Y_{1})(u-1)\}(u-1)\Big].
    \end{align*}

      Here is the main assertion for the test statistic $S_T$
      under the null hypothesis ${\cal{H}} _{02} $:

      \begin{thm} 
  \label{thm_2_4}Let $\{Y_t\}_{t\in \mathbb N}$ follow   the model~\eqref{inarch} satisfying the null hypothesis ${\cal{H}}_{02}$. Assume that the estimators $\widehat \theta_{1T},\widehat \theta_{2T}$ satisfy (\ref{INARCHest}),  and that
  \begin{equation}\label{ass}
  \mathbb {E} \left( \int_0^1 \Big[1+\exp\{2( \theta_1 +\theta_2 Y_1)(u-1)\}(1+Y^4_{1})\Big]w(u)du\right)<\infty.
    \end{equation}

  Then
  the limit distribution of $S_T$ defined in \eqref{testINAR} is the same as that of
   \[
       \int_0^1 \bar{\cal{V}}_T^2 (u;
       \theta_{10}, \theta_{20}) w(u) du,
       \]
       where
$\bar{\cal {V}}_T(u;\theta_{10},\theta_{20})=\frac{1}{\sqrt T}\sum_{t=q+1}^T  {\cal V}_t(u;\theta_{10},\theta_{20})$
with $ \{{\cal V}_t(u;\theta_1,\theta_2), \, u\in [0,1]\}$ as in eqn.~\eqref{calVt}. 
 Moreover, the process $\{ \bar{\cal{V}}_T (u;\theta_{10},\theta_{20});\ u\in[0,1]\}$  converges in distribution to a zero--mean Gaussian process
    $\{{\cal{V}} (u;\theta_{10},\theta_{20});\ u\in[0,1]\}$
    with covariance structure $\mathbb {E}\Big[{\cal{V}}_{q+1}(u_1;\theta_{10},\theta_{20}){\cal{V}}_{q+1}(u_2;\theta_{10},\theta_{20})\Big]$, $u_1, u_2 \in [0,1]$
and $S_T$ converges in distribution to
  \[
      \int_0^1 {\cal{V}}^2 (u;
       \theta_{10}, \theta_{20}) w(u) du.
       \]
      \end{thm}
The proof is postponed to the Appendix.

\section{Estimation of parameters and bootstrap test}\label{sec5}
Recall that the test statistic $S_T$ suggested in Section~\ref{sec_1} implicitly depends on estimated parameters, and that the asymptotic null distribution of $S_T$ derived in Section~\ref{sec3} assumes certain properties for these estimators (see (A.4) for the INAR model and \eqref{INARCHest} for the INARCH model). There  is a number of estimation methods with corresponding estimators having the desired properties. Here we construct estimators of the parameters based on conditional least squares  (CLS) along the line of \cite{klimko}. Again the focus is on the INAR(1) model since  the estimators for the INARCH(1) Poisson
 model are introduced analogously and the computational expressions are also similar.

  To begin with notice that under the INAR(1) \eqref{modelINAR} satisfying  ${\cal{H}}_{01}$ with  parameter $\vartheta=(p, \theta)'$ we have
\begin{equation}\label{eq51}
\mathbb{E}_{\vartheta}( Y_t| Y_{t-1})= pY_{t-1}  +
\mathbb{E}_{\theta}(\varepsilon_t).
\end{equation}
The CLS estimator  $\widehat\vartheta_T=(\widehat p_T, \widehat \theta_T)'$  of $\vartheta=(p,\theta)'$  is defined as a minimizer of
  \begin{equation}\label{eq52}
  L(\vartheta)=\sum_{t=2}^T\Big(Y_t-  pY_{t-1}  -
\mathbb{E}_{\theta}(\varepsilon_t) \Big)^2
  \end{equation}
 w.r.t. $\vartheta$.

Analogously, the CLS estimator in the INARCH(1) Poisson model  satisfying
 ${\cal{H}}_{02}$ is defined as in \eqref{eq52}, but using the equation
\[
\mathbb{E}_{\vartheta}( Y_t| Y_{t-1})=\theta_1+\theta_2 Y_{t-1}
\]
instead of \eqref{eq51}.

   Concerning limit properties of  CLS  estimator in INAR(1)  under ${\cal{H}}_{01}$ and  assumptions (A.1) -- (A.3), as
$T \to \infty,$
 \begin{align*}
    \sqrt T &(\widehat p_T-p_0, \widehat \theta_T-\theta_0)^{'}\\
&=\frac{1}{\sqrt T}\boldsymbol B^{-1}(\vartheta_0)
\Big(\sum_{t=2}^T\Big(Y_t-\mathbb E_{\vartheta_0}(Y_t|Y_{t-1})\Big) \Big(Y_{t-1}, \frac{\partial \mathbb{E}_{\theta}
\varepsilon_t }{\partial \theta}\Big|_{\theta=\theta_0}\Big)^{'}
+o_P(1),
    \end{align*}
    where
    $$
    \boldsymbol B(\vartheta)= \begin{pmatrix}   \mathbb{E} Y^2_1 &  \mathbb{E}Y_{t-1}\frac{\partial \mathbb{E}_{\theta} \varepsilon_t }{\partial \theta}\\
  \mathbb{E} Y_{t-1}\frac{\partial \mathbb{E}_{\theta} \varepsilon_t }{\partial \theta} &
    \Big(\frac{\partial \mathbb{E}_{\theta} \varepsilon_t }{\partial \theta}\Big)^2
\end{pmatrix},
    $$
which immediately implies that the CLS estimators have the property  (A.4). The derivation follows closely lines  of  those in \cite{klimko} and therefore, we omit details.

 Now we shortly discuss behavior  of the  CLS estimators $(\widehat p_T,\widehat \theta_T)'$
under alternatives.  Recall  $(\widehat p_T,\widehat \theta_T)'$  are minimizers of (\ref{eq51})
where  $\mathbb{E}_{\theta}\varepsilon _t$  is  the expectation under the null hypothesis, however under alternatives we  have  generally $\mathbb{E}\varepsilon _t$. Denoting $p_0$ the true parameter value  we have a look at minimizers of
\begin{align*}
\mathbb{E} &\left[L(p, \theta)\right]\frac{1}{T}=
\Big\{\mathbb E(Y_t- E (Y_t|Y_{t-1}))^2  + \mathbb{E} ( p_0Y_{t-1} +\mathbb{E}\varepsilon _t -pY_{t-1} - \mathbb {E}\varepsilon _t)^2 (1+o(1))\Big\}\\
&= \Big\{\mathbb {E}(Y_t- E (Y_t|Y_{t-1}))^2 +(p_0-p)^2\mathbb{E} (Y_{t-1})^2\Big\}+(\mathbb{E}_{\theta}\varepsilon _t -\mathbb{E}\varepsilon _t )^2.
\end{align*}
 It is easily seen that minimum is reached for $p=p_0$  and for $\theta_A$ that minimizes
  $(E\varepsilon _t - E_{\theta}\varepsilon _t)^2$ w.r.t. $\theta$. If such $\theta_A$ exists we get parallel to the case of the null hypothesis that
$$
  \widehat p_T-p=o_P(1)
  , \quad
  \widehat \theta_T-\theta_A=o_P(1).
  $$
   Hence  if there exists $\theta_A$  minimizing
  $(E\varepsilon _t - E_{\theta}\varepsilon _t)^2$ w.r.t. $\theta$, then  the   CLS estimators $(\widehat p_T,\widehat \theta_T)$ have the property required in Theorem  4.3.

\medskip

As it was shown  in Section~\ref{sec3},  the asymptotic null distribution of $S_T$ is complicated  and depends on several unknown quantities including the true value of the parameter $\vartheta$. Therefore, some resampling scheme is adopted in order to carry out the test procedure and compute critical points. In what follows we advocate the parametric bootstrap as resampling scheme because it reflects all aspects of the underlying model, and has been put on a firm theoretical basis both with i.i.d. data, \cite{genest}, as well as with data involving dependence, \cite{leucht}.

We shall outline the parametric bootstrap for the INAR model, the corresponding procedure for the INARCH model being completely analogous. Specifically in view of the data $Y_t, \ t=1,...,T$, and in order to carry out the test we compute the parameter estimate $\widehat \vartheta_T=(\widehat p_T,\widehat \theta_T)'$, and the corresponding value of the test statistic $S_T:=S_T(Y_1,...,Y_T;\widehat \vartheta_T)$. Then the parametric bootstrap takes the following form:

\begin{enumerate}
\item \label{enum:1_1} Generate realizations $U^*_{t,i}, \ i=1,...$, where $U^*_{t,i}$ are as $U_{t,i}$ in eqn.~\eqref{modelINAR} but with $p$ replaced by  $\widehat p_T$.
\item \label{enum:1_2} Generate realizations $\varepsilon^*_t, \ t=1,...$, where $\varepsilon^*_t$ are as $\varepsilon_t$ in eqn.~\eqref{modelINAR} but with $\theta$ replaced by  $\widehat \theta_T$.
\item \label{enum:1_3} Compute pseudo--observations $Y^*_t, \ t=1,...$, using eqn.~\eqref{modelINAR} with $(U_{t,i},\varepsilon_t)$ replaced by $(U^*_{t,i},\varepsilon^*_t)$.
\item \label{enum:2}  Fit the model \eqref{modelINAR} again with ${{Y}}_t^*, \ t=1,...,T$, as observations to obtain the estimate ${\widehat\vartheta}_T^*=({\widehat{ p}}_T^*,{\widehat\theta}_T^*)'$.
\item \label{enum:3} Compute the test statistic $S^*_T= S_T(Y^*_1,...,Y^*_T;\widehat \vartheta_T^*)$.
\item Steps (1) to (5) are then repeated $B$ times to obtain the sequence of test statistics, say, $S^*_{1,T}, \dots, S^*_{B,T}$.
\end{enumerate}
Let $ S^*_{(1)}\leq S^*_{(2)}\leq \cdots \leq S^*_{(B)}$, be the corresponding order statistics. Then the null hypothesis is rejected at level of significance $\alpha$ if the value of the test statistic based on the original data exceeds the $(1-\alpha)B$ empirical quantile obtained from $S^*_{1,T}, \dots, S^*_{B,T}$, i.e. when $S_T>S^*_{(B-\alpha B)}$.

\section{Extension to higher order}

We discuss possible extension of the procedures to higher order. By way of example, we  consider the INAR(2) model formally defined by the equation
\[ Y_t=\sum_{i=1}^{Y_{t-1}}U_{t,1i}+\sum_{i=1}^{Y_{t-2}}U_{t,2i}
+\varepsilon_t,\quad t=1,\ldots,
\]
where for $i=1,2$, \  $U_{t,im}, m=1,2,...,$ are i.i.d.  random variables  with finite variance  such that
$P(U_{t,im}=1)=p_i,\, P(U_{t,im}=0)=1-p_i,$ $p_i\in[0,1]$, $p_1+p_2<1$ and $U_{im}$ are independent of $Y_{t-i}$,
the sequences $U_{t,1m}$ and $U_{t,2m}, \ m=1,2,...$, are mutually independent, and the i.i.d. innovations $\varepsilon_t,\, t=1,\ldots ,$ have finite second moment and are independent of
  $Y_{t-i}, i=1,2$.  Following the lead of eqn. (7), we have
\begin{align} \label{PGF2}
 g_{Y_t}(u)=&\mathbb E\left[u^{Y_t}\right]=\mathbb E \left [\mathbb E(u^{Y_t}| Y_{t-1},Y_{t-2})\right]\\
= \nonumber & g_{\varepsilon}(u; \theta) g_{Y_{t-1},Y_{t-2}}( 1+p_1(u-1),1+p_2(u-1)), \quad u\in [0,1],
\end{align}
 where $g_{Y_{t},Y_{t-1}}(u,v)$ denotes the joint PGF of $Y_t$ and $Y_{t-1}$. Under the above assumptions $\{Y_t\}_{t\in\mathbb{N}}$ is stationary and $g_{Y_{t},Y_{t-1}}$ may be estimated by the (joint) empirical PGF
\begin{equation} \label{EPGF2}
\widehat g_{T,Y}(u,v)=\frac{1}{T-1} \sum_{t=2}^T u^{Y_{t}}v^{Y_{t-1}}.
\end{equation}  

 Based on eqn. (\ref{PGF2})  and given suitable estimators $\widehat \vartheta_T=(\widehat p_{1T},\widehat p_{2T},\widehat\theta_T)'$   a natural semiparametric estimate of the joint PGF is
\begin{equation} \label{EPGF2null}
\widehat g_{T0}(u)=g_{\varepsilon}(u,\widehat \theta_T)   \widehat
g_{T,Y}(1+\widehat p_{1T}(u-1), 1+\widehat p_{2T}(u-1)), \quad u\in \mathbb{R},
\end{equation}
where $\widehat g_{T,Y}(1+\widehat p_{1T}(u-1), 1+\widehat p_{2T}(u-1))$  is the empirical PGF in eqn. (\ref{EPGF2}) computed at the point $(u,v)=(1+\widehat p_{1T}(u-1),1+\widehat p_{2T}(u-1))$.
Hence a test statistic analogous to that of eqn.~\eqref{testINAR} may be defined by using the quantities in eqns. (\ref{EPGF2}) and (\ref{EPGF2null}), instead of those in eqns. (6) and (8), respectively.  The case of the INARCH(2) and that of higher order models can be treated analogously. Moreover, based on assumptions analogous to those of Section 4, the asymptotic results also go through on the grounds of entirely parallel reasoning. Therefore we do not pursue this here in more detail in order to save space. As a last note, clearly  computations become somewhat cumbersome with increasing model order, but in principle calculating the test statistic even by numerical integration should not be a problem.

\section{Simulations}

In this section we study the small--sample behavior of the suggested bootstrap test via a simulation study.

\subsection{{\rm{{\bf{INAR(1) and INARCH(1)}}}}}

We consider the null hypotheses of Poisson INAR(1) and Poisson INARCH(1) models and investigate the size of the test under the null hypothesis as well as the power under various alternatives. The test statistic $S_T$ is computed using CLS estimators of the model parameters. The  weight function $w(u)$ used is $w(u)=u^a$ for $a=0,1,2,5$. The  p-value  of the bootstrap test is computed from  $B=499$ bootstrap samples and the percentage of rejection of the test  is estimated from $500$ repetitions. The simulations were conducted in the R-computing environment \cite{rko}. In the following we present only a part of the results for some particular cases.   Results for additional settings (leading to mostly analogous conclusions) could be provided by the authors upon a request.

\begin{figure}
\includegraphics[width=0.48\textwidth]{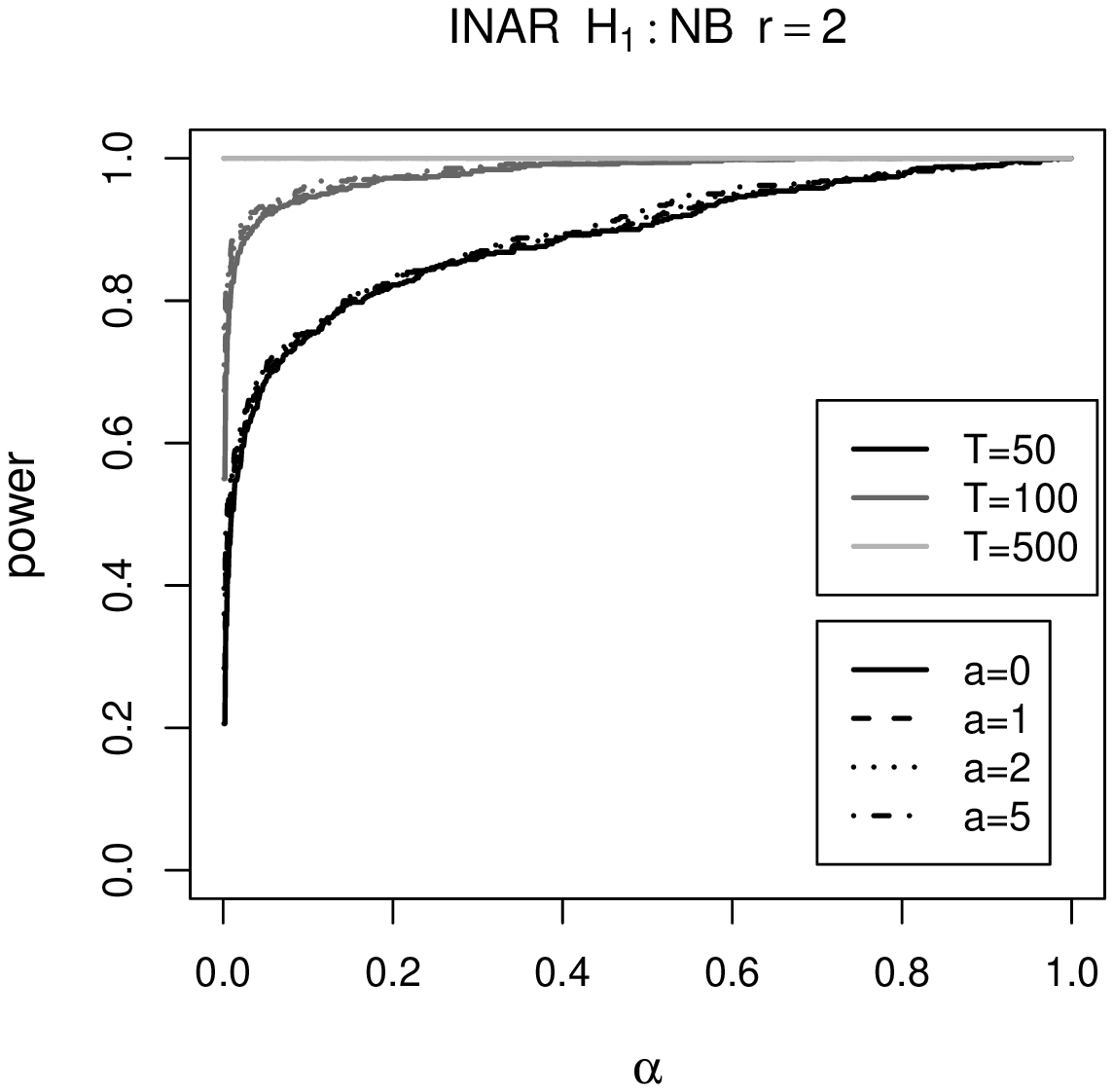}
\includegraphics[width=0.48\textwidth]{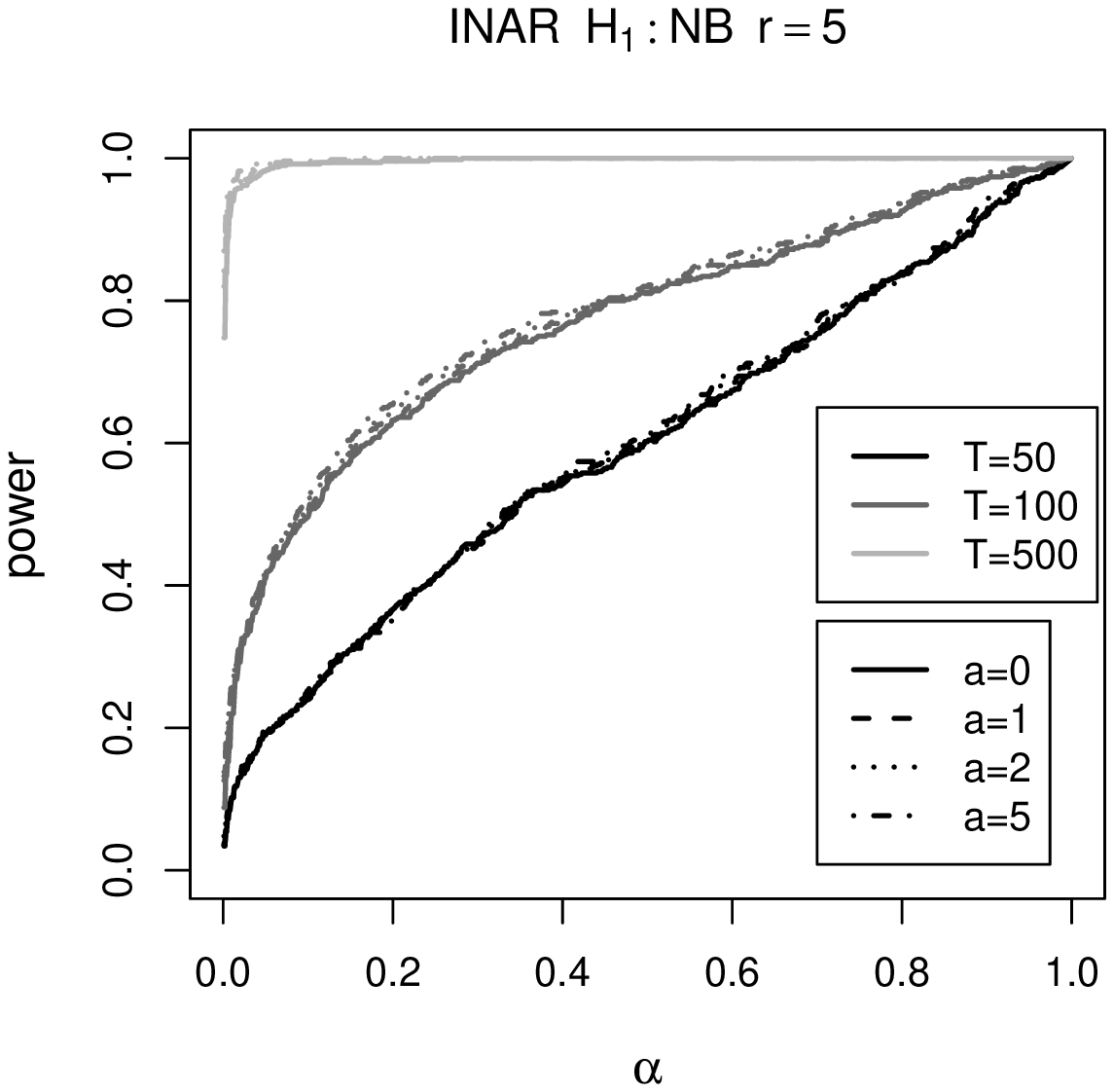}
\caption{Power of the test for a Poisson INAR(1) against an INAR(1) model with Negative Binomial innovations with dispersion parameter $r$. }\label{fig-inar-nb}
\end{figure}

\begin{figure}
\includegraphics[width=0.48\textwidth]{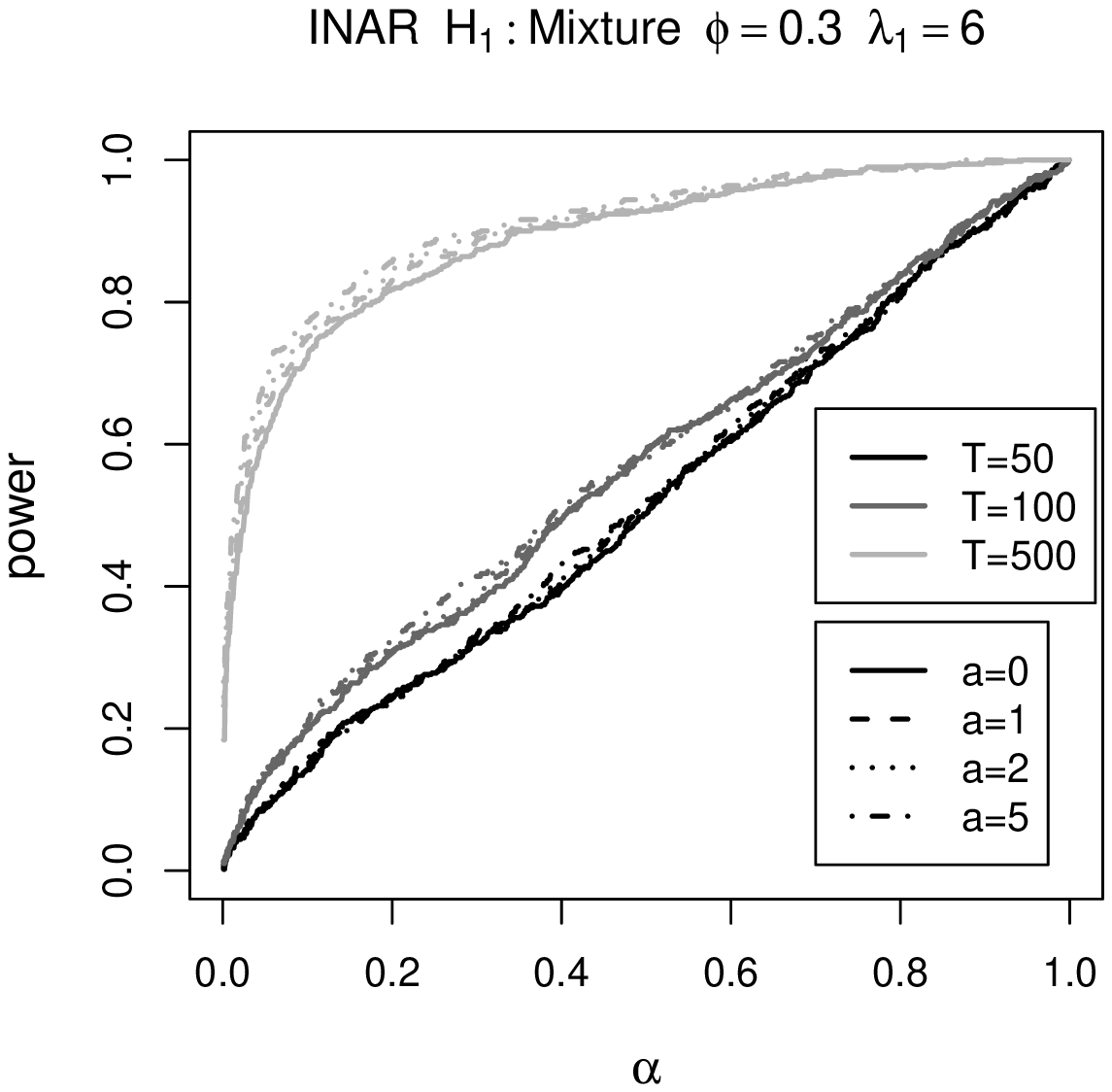}
\includegraphics[width=0.48\textwidth]{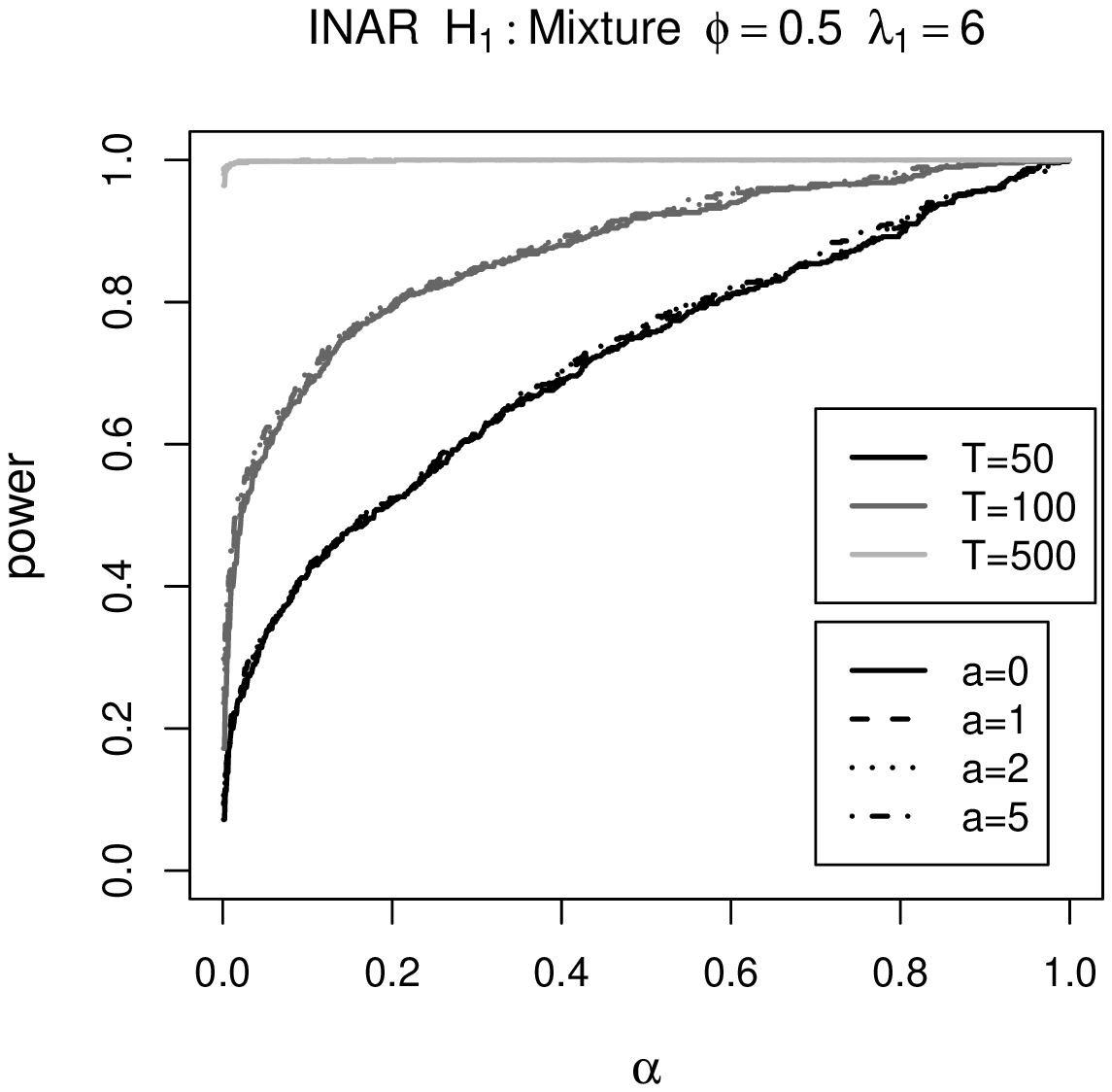}
\caption{Power of the test for a Poisson INAR(1)  under  an INAR(1) with innovations following a mixture of two Poisson distributions with mixing parameter $\phi$. }\label{fig-inar-mix}
\end{figure}

\begin{figure}
\includegraphics[width=0.48\textwidth]{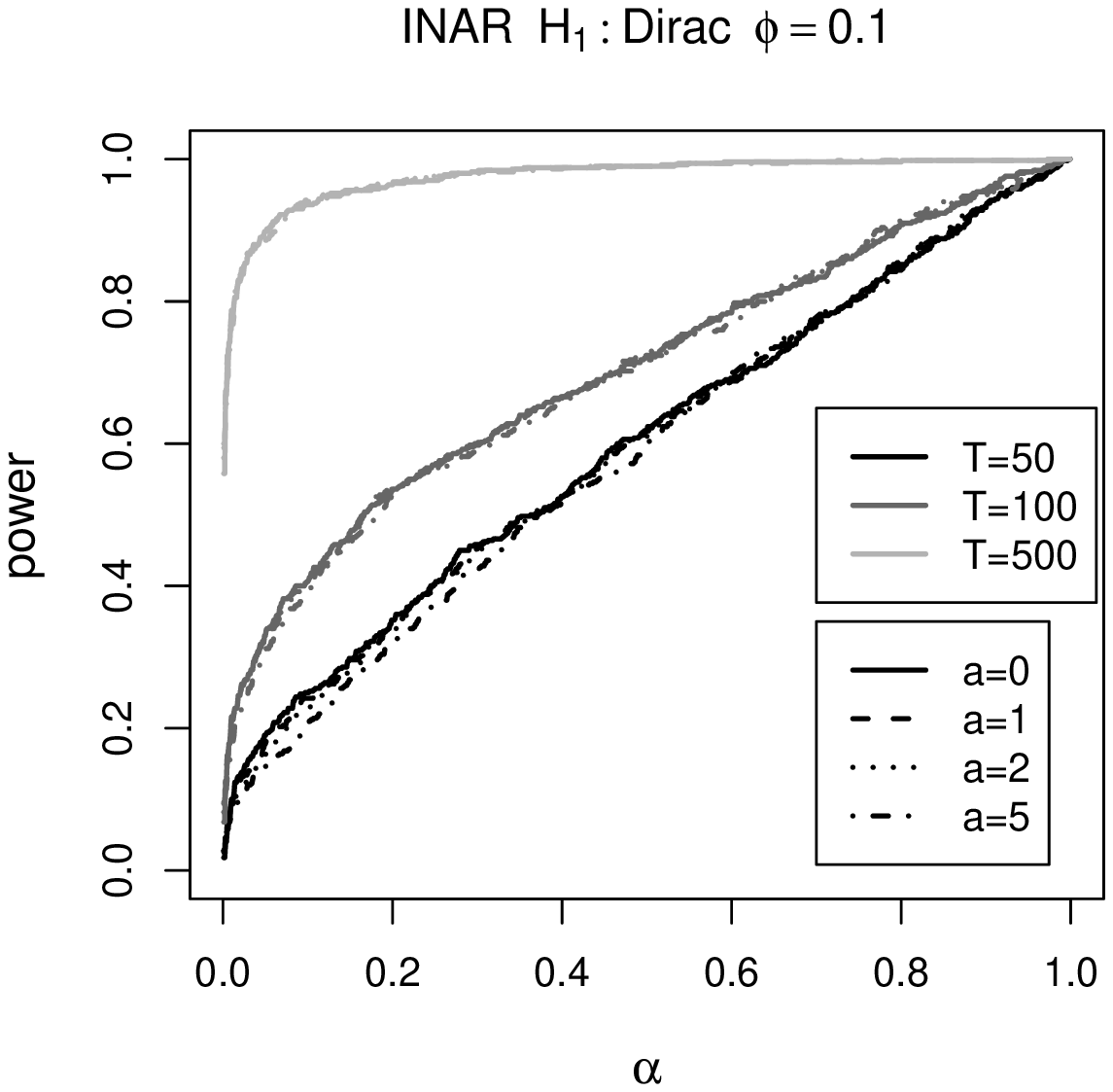}
\includegraphics[width=0.48\textwidth]{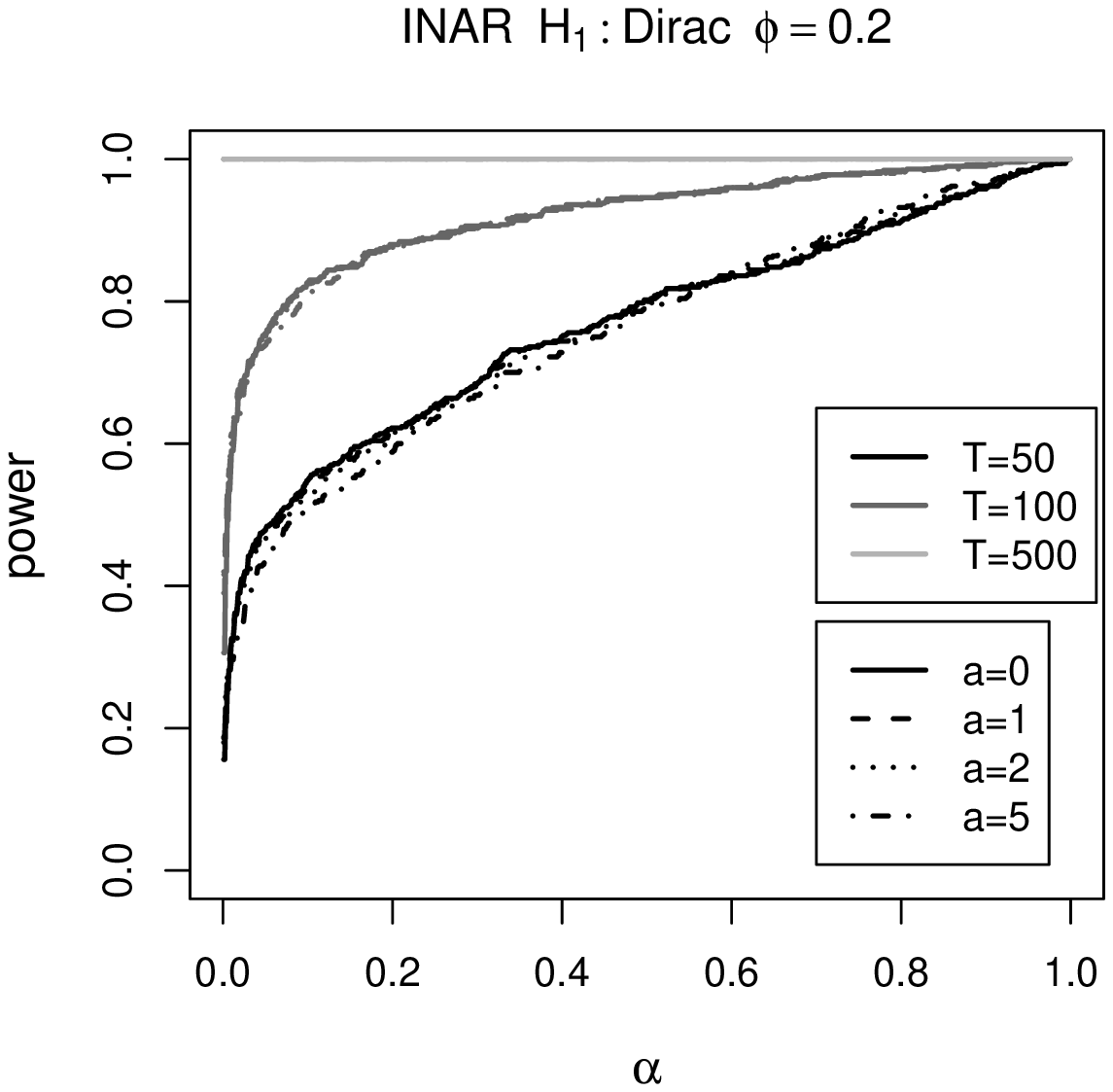}
  \caption{Power of the   test for a Poisson INAR(1)   under an INAR(1) with innovations following a mixture of a Poisson distribution, and a Dirac measure at 0 with weight $\phi$.}  
 \label{fig-inar-dirac}
\end{figure}

\begin{figure}
\includegraphics[width=0.48\textwidth]{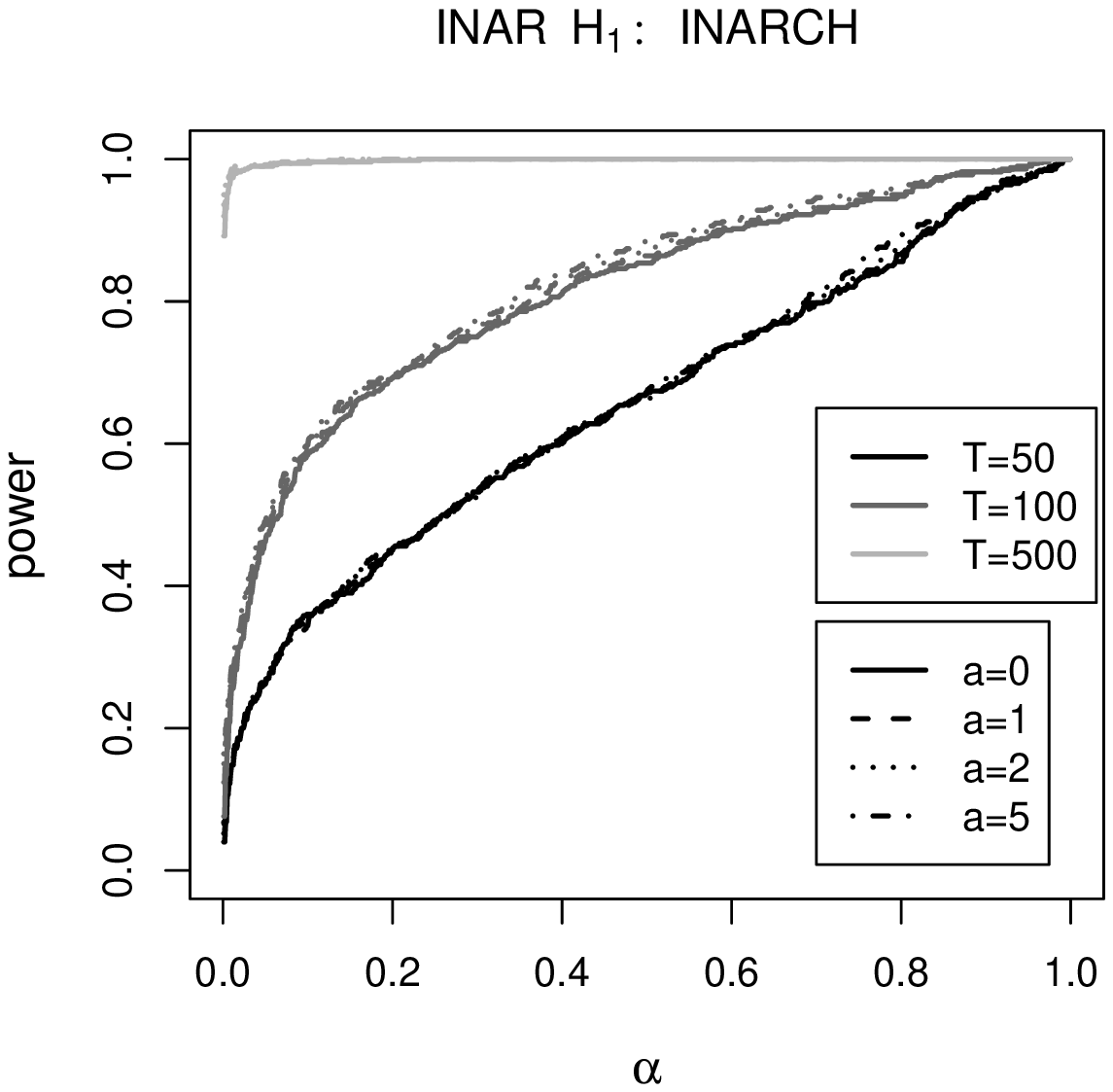}
\includegraphics[width=0.48\textwidth]{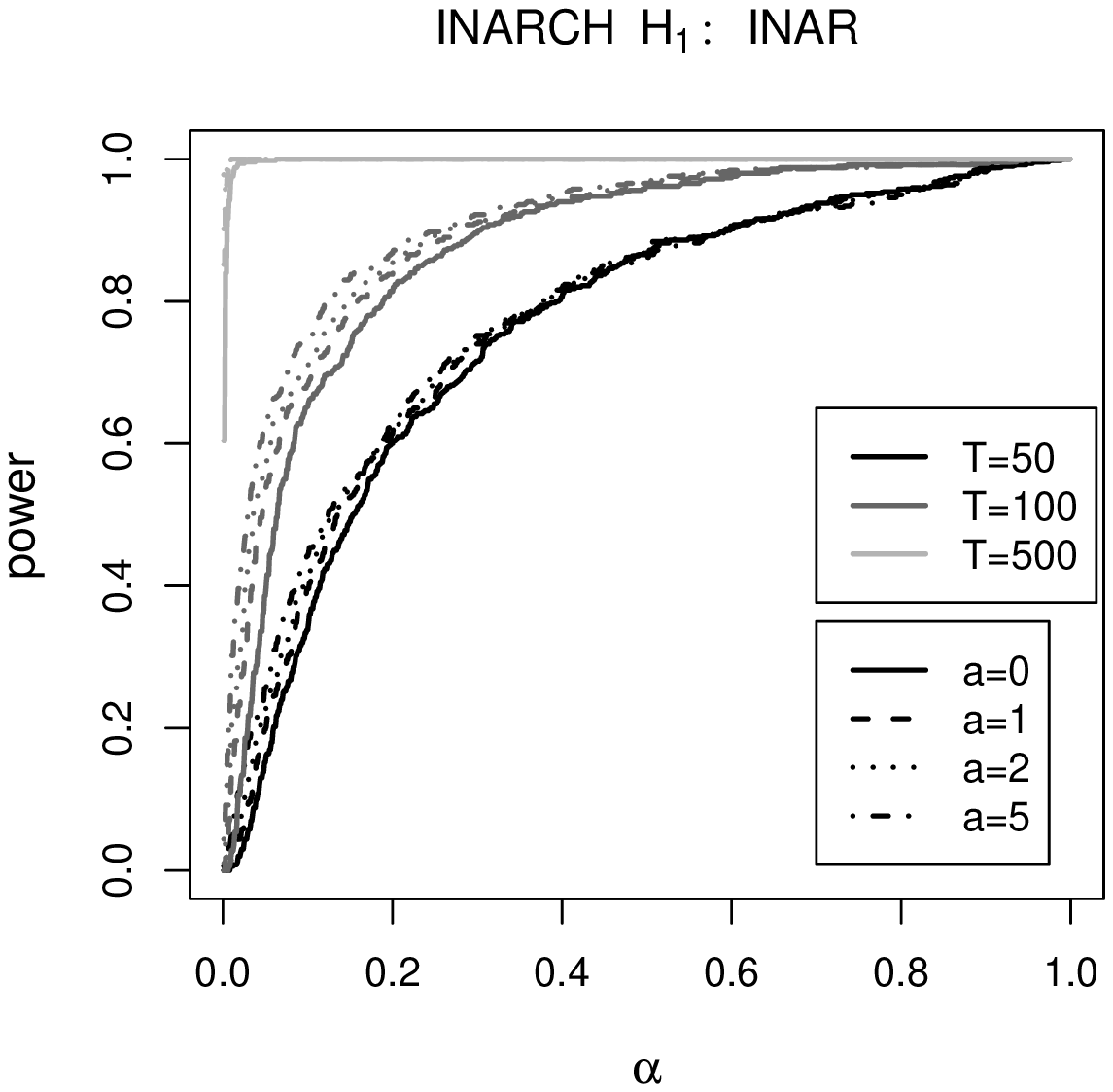}
 \caption{Power of the  test for a  Poisson INAR(1) under a Poisson INARCH(1)  (left panel). Power for a Poisson INARCH(1) under a Poisson INAR(1) (right panel). }\label{fig-inar-inarch}
\end{figure}

\begin{figure}
\includegraphics[width=0.48\textwidth]{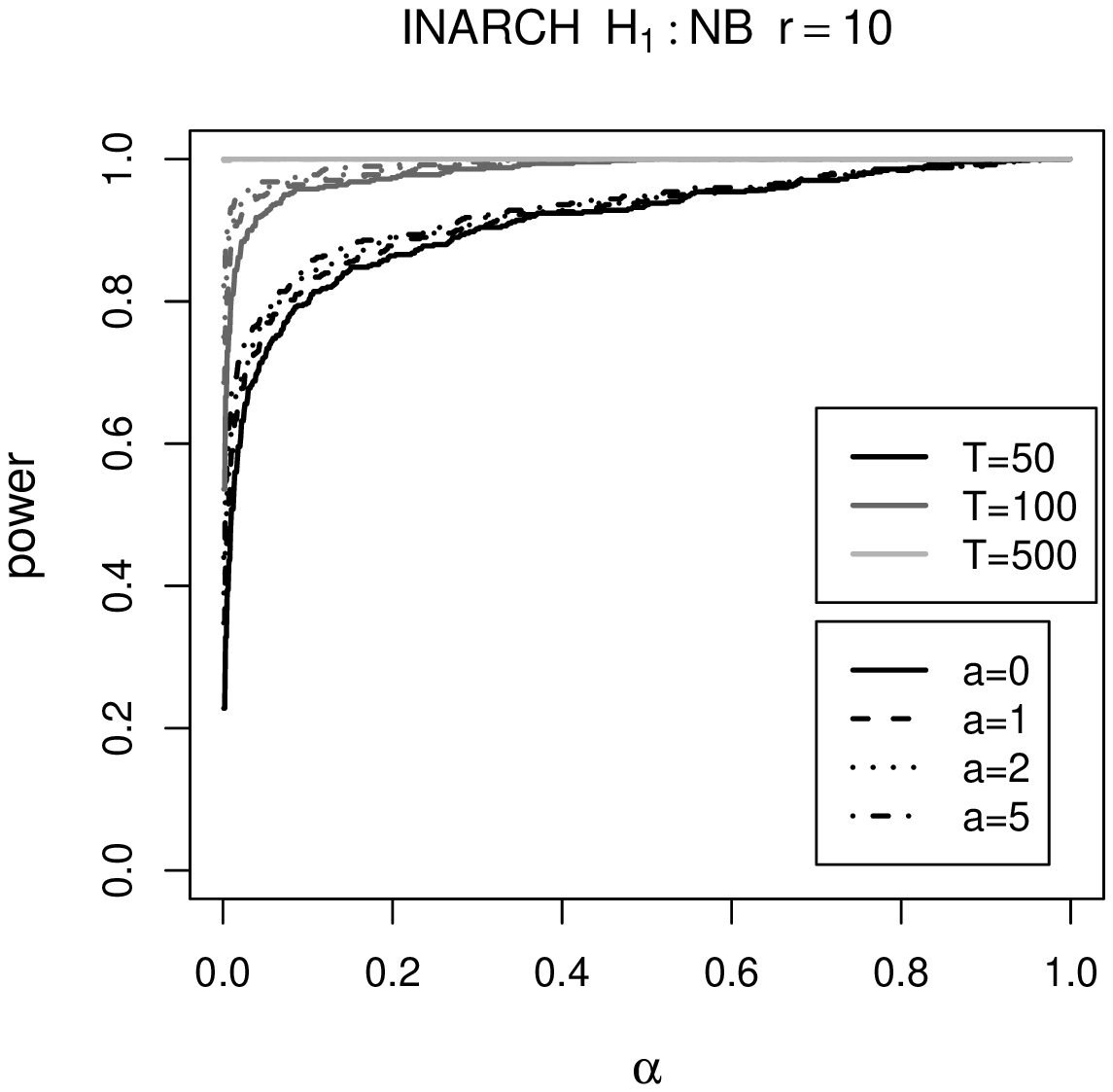}
\includegraphics[width=0.48\textwidth]{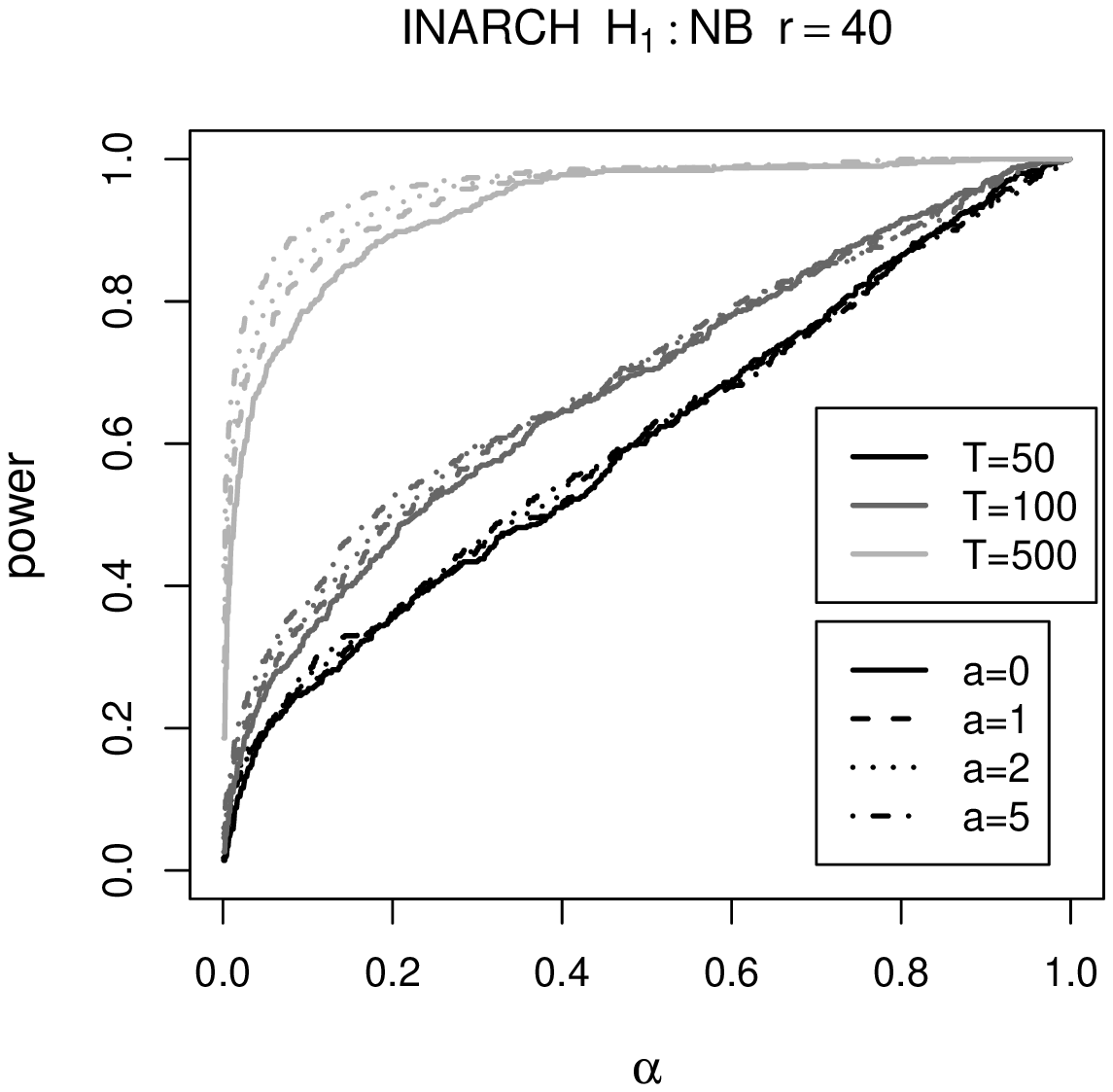}
\caption{Power of the  test for a Poisson INARCH(1)  under a Negative Binomial INARCH(1)  with  dispersion parameter $r=10$ and $r=40$. }\label{fig-inarch-nb}
\end{figure}

\begin{figure}
\centering
\includegraphics[width=0.48\textwidth]{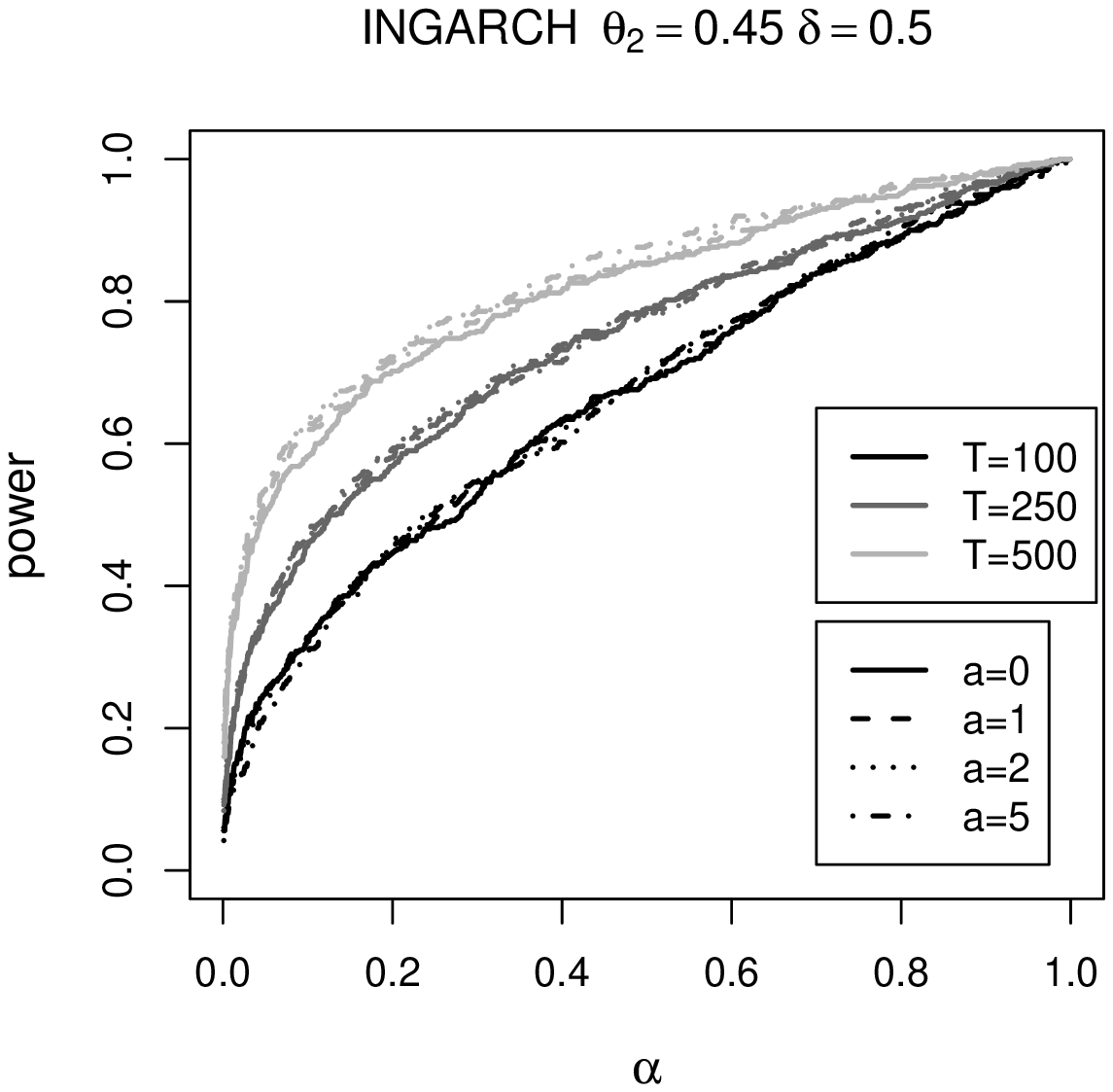}
\includegraphics[width=0.48\textwidth]{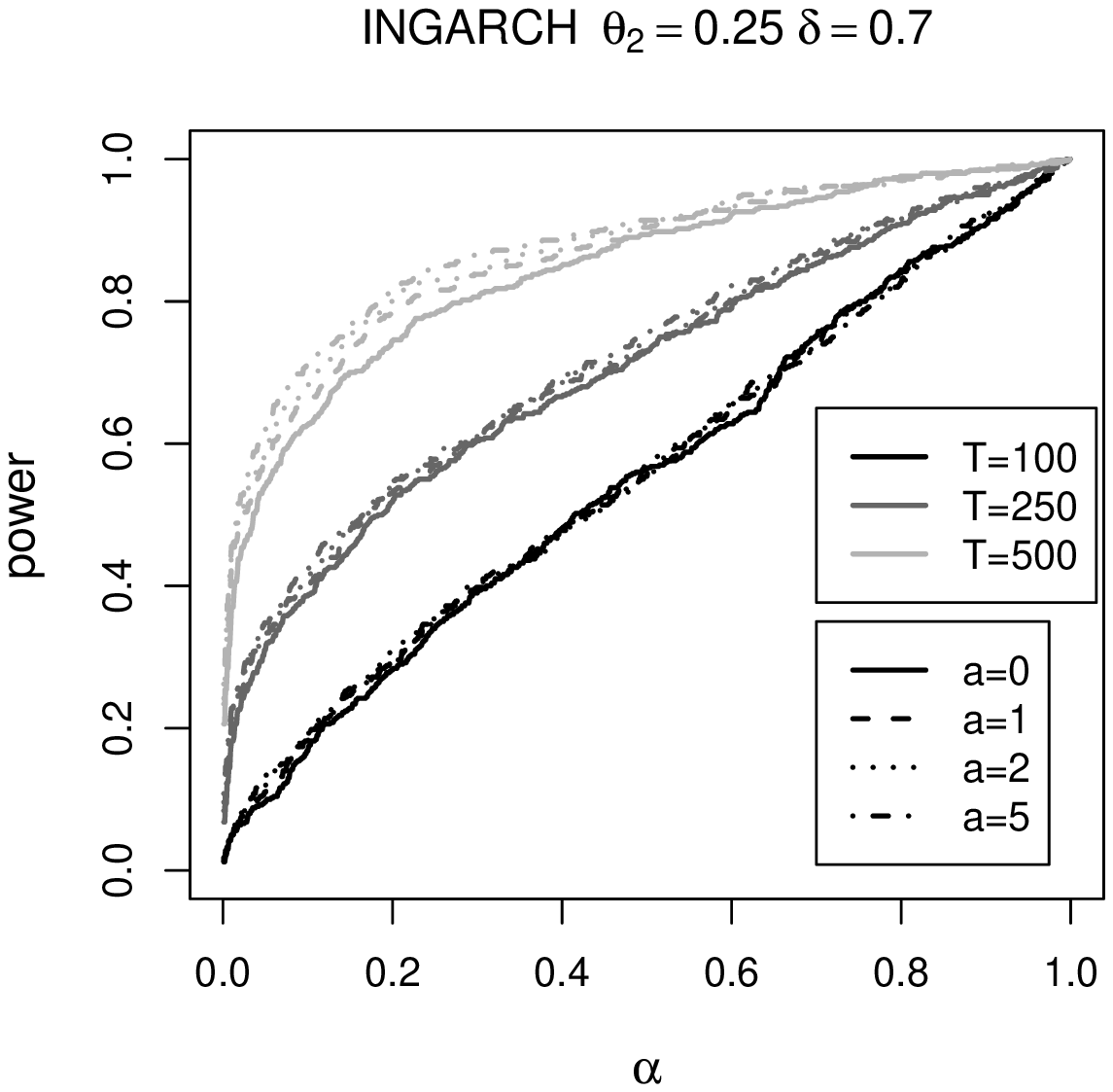}
\caption{ Power of the test for a Poisson INARCH(1) under a Poisson INGARCH(1,1) alternative with parameters $\theta_2$ and $\delta$.}\label{fig-inarch-garch}
\end{figure}

\begin{figure}
\centering
\includegraphics[width=0.48\textwidth]{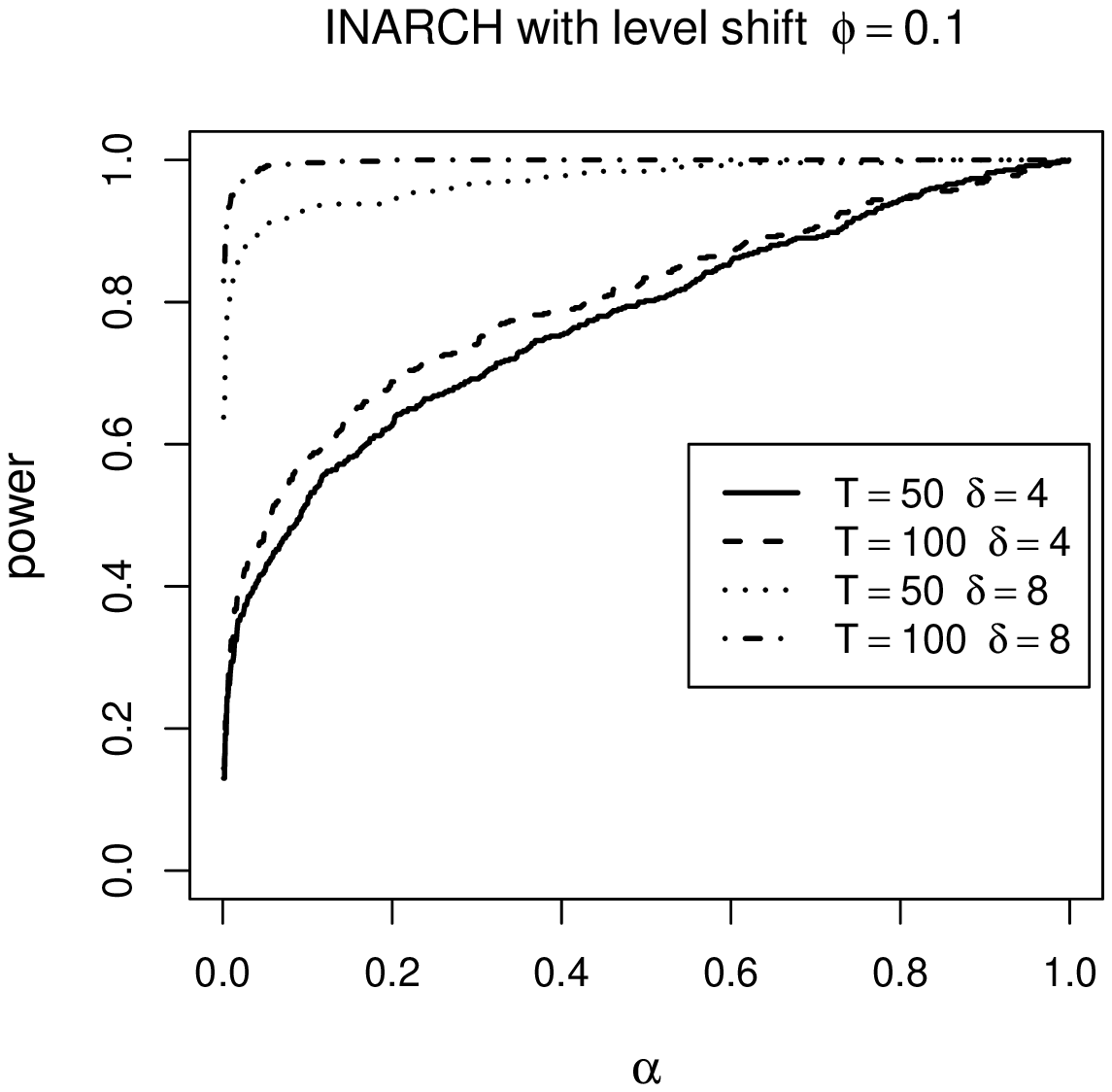}
\includegraphics[width=0.48\textwidth]{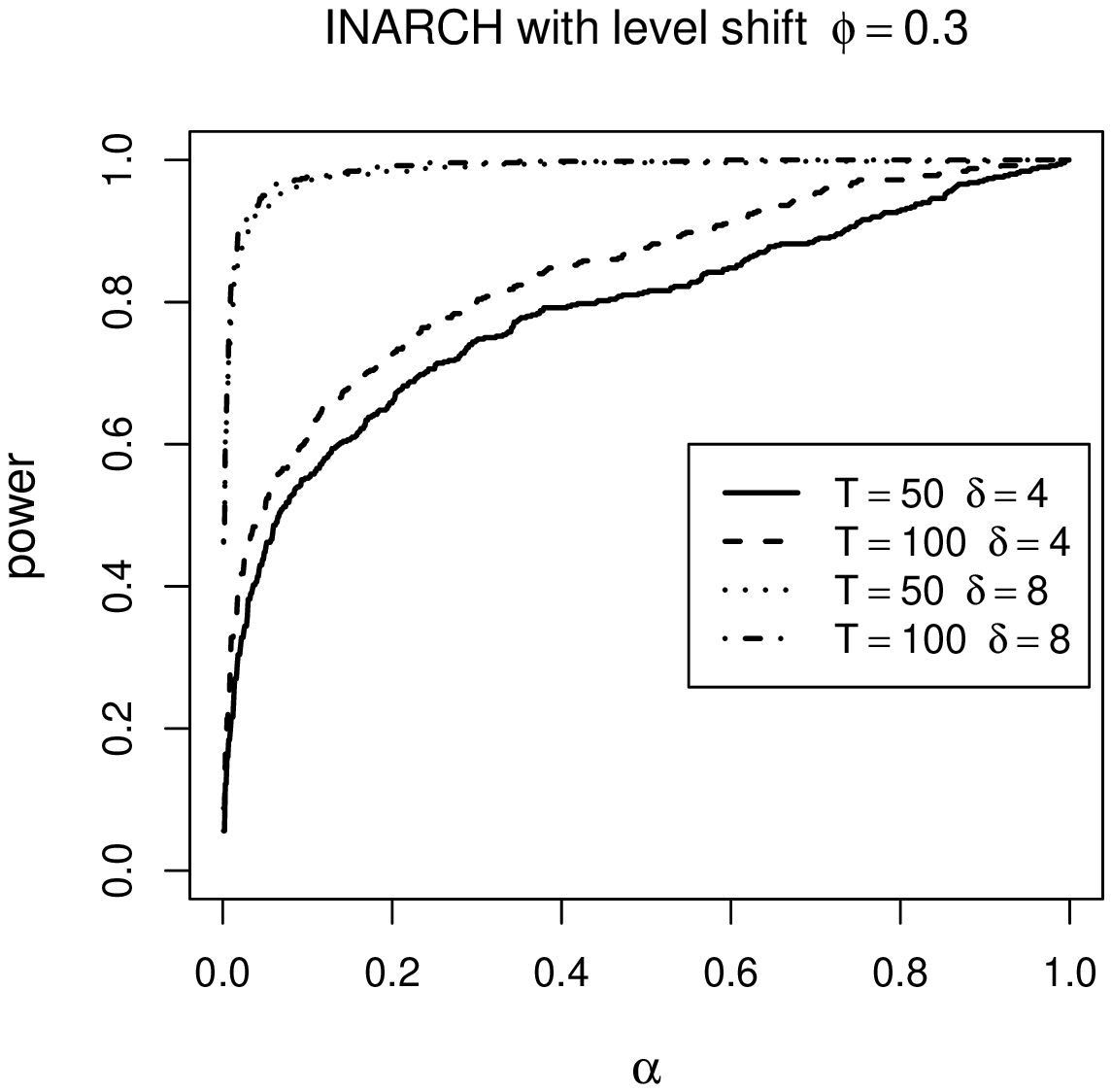}
\caption{ Power of the test for a Poisson INARCH(1) under a  Poisson INARCH(1) with level shift  with parameters $\delta=4$ or $\delta=8$, and $\phi=0.1$ or $\phi=0.3$.}\label{fig-inarch-ls}
\end{figure}

Specifically, we present the observed percentage of rejection under the null hypotheses ${\cal{H}}_{01}$ and ${\cal{H}}_{02}$. In the former case we use an  INAR(1) model with $p=0.6$ under Poisson innovations $\eps_t$ with mean $\theta=4$, throughout. In turn, under the null hypothesis ${\cal{H}}_{02}$ a Poisson INARCH(1) model with $\theta_1=4$ and $\theta_2=0.6$ was employed in all cases. The Monte Carlo sample sizes used are $T=50, 100, 250$, and $500$. In Table~\ref{t1}, the realized size of both tests is shown corresponding to significance levels $\alpha=0.01, 0.05$, and $0.1$. Power results  for the test for the null hypothesis ${\cal{H}}_{01}$ (resp. ${\cal{H}}_{02}$) for various alternatives are shown in Figures \ref{fig-inar-nb}--\ref{fig-inar-dirac} (resp. \ref{fig-inarch-nb}--\ref{fig-inarch-ls}), while Figure \ref{fig-inar-inarch} shows the power of the test for the null hypothesis ${\cal{H}}_{01}$ against an INARCH model, and the power of the test for the null hypothesis ${\cal{H}}_{02}$  against an INAR model, in all cases with the parameter values just mentioned.

\begin{table}[htb]
\centering
\tbl{Size of the test for the Poisson INAR(1) model (left) and the test for the Poisson INARCH(1) model (right).}
{\begin{tabular}{rr|rrr}
\toprule
  \multicolumn{5}{c}{$H_0:$ INAR(1)}\\
  \colrule
&&\multicolumn{3}{c}{$\alpha$}\\
$T$ & $a$ & 0.01 & 0.05 & 0.1 \\
\colrule
50 & 0 & 0.008 & 0.044 & 0.100 \\
  50 & 1 & 0.008 & 0.050 & 0.102 \\
  50 & 2 & 0.006 & 0.054 & 0.112 \\
  50 & 5 & 0.008 & 0.056 & 0.110 \\
  100 & 0 & 0.004 & 0.038 & 0.080 \\
  100 & 1 & 0.004 & 0.040 & 0.082 \\
  100 & 2 & 0.006 & 0.042 & 0.088 \\
  100 & 5 & 0.012 & 0.052 & 0.102 \\
  500 & 0 & 0.018 & 0.038 & 0.090 \\
  500 & 1 & 0.016 & 0.044 & 0.090 \\
  500 & 2 & 0.010 & 0.042 & 0.088 \\
  500 & 5 & 0.010 & 0.044 & 0.094 \\
   \botrule
\end{tabular}
\hskip5mm
\begin{tabular}{rr|rrr}
  \toprule
    \multicolumn{5}{c}{$H_0:$ INARCH(1)}\\
  \colrule
&&\multicolumn{3}{c}{$\alpha$}\\
$T$ & $a$ & 0.01 & 0.05 & 0.1 \\
  \colrule
50 & 0 & 0.006 & 0.038 & 0.082 \\
  50 & 1 & 0.008 & 0.036 & 0.092 \\
  50 & 2 & 0.006 & 0.040 & 0.090 \\
  50 & 5 & 0.004 & 0.044 & 0.090 \\
  100 & 0 & 0.002 & 0.036 & 0.094 \\
  100 & 1 & 0.004 & 0.042 & 0.096 \\
  100 & 2 & 0.002 & 0.042 & 0.096 \\
  100 & 5 & 0.006 & 0.044 & 0.094 \\
  500 & 0 & 0.006 & 0.036 & 0.080 \\
  500 & 1 & 0.006 & 0.034 & 0.092 \\
  500 & 2 & 0.004 & 0.042 & 0.092 \\
  500 & 5 & 0.008 & 0.040 & 0.102 \\
   \botrule
\end{tabular}}
\label{t1}
\end{table}


For the null hypothesis of Poisson INAR(1) we considered four different sets of possible alternatives. Namely, an INAR(1) with the innovations $\eps_t$ following

(a) a  Negative Binomial distribution,

(b) a mixture of two Poisson distributions, and

 (c) a mixture of a Poisson and a Dirac measure at $0$,

\noindent all with mean $\theta=4$, and

 (d) the Poisson INARCH(1),

\noindent as a model--deviation alternative.

  All these models serve as possible alternatives to the Poisson INAR(1) model for data which exhibit an  overdispersion.
 
 The results in Table~\ref{t1} suggest that the bootstrap test, despite being mildly under--sized or over--sized, it generally keeps the prescribed significance level to a satisfactory degree.

The power for the alternative (a) is plotted as a function of the significance level $\alpha$ in Figure~\ref{fig-inar-nb}.
Here, the innovations $\eps_t$ are generated from a Negative Binomial distribution with a dispersion parameter $r=2$ and $r=5$ respectively (i.e. ${\rm{Var}}( \eps_t)=\theta(1+\theta/r)$). For $r=2$ we obtain a reasonable power already for small sample size ($T=50$). However as $r$ growths, the innovation distribution tends to the Poisson distribution and one needs to have more observation in order to obtain sufficient power. As an example at level $\alpha=0.05$, for $r=5$ and sample size $T=100$ we observe power only around 40~\%, but $T=500$ leads to percentage of rejection close to 100~\%.

Under alternative (b), the innovations are generated from  a mixture of two Poisson distributions (always with mixture mean equal to 4) of the form $\phi\mathsf{Po}(\lambda_1)+(1-\phi)\mathsf{Po}(\lambda_2)$, where $\mathsf{Po}(\lambda)$ denotes the Poisson distribution with mean $\lambda$.  The estimated power for $\lambda_1=6$ and $\phi=0.3$ or $\phi=0.5$ is plotted in~Figure~\ref{fig-inar-mix}.
We can see that  $\phi=0.5$ leads to noticeably larger power compared to $\phi=0.3$.
Likewise for a larger $\lambda_1$ one would obtain a larger power (results not shown). On the other hand, if $|\theta-\lambda_1|$ decreases then the mixture distribution gets closer to the Poisson distribution  and consequently a relatively long series is needed in order to obtain a reasonable power.

The alternative (c) is plotted in Figure~\ref{fig-inar-dirac}. Here,  the innovations are assumed to follow a mixture $\phi\mathsf{D}_0+(1-\phi)\mathsf{Po}(\lambda)$, where $\mathsf{D}_0$ denotes the
Dirac measure at $0$. As expected, the power increases with the value of $\phi$.

 Finally, Figure~\ref{fig-inar-inarch} (left panel) shows  that the test is also able to distinguish between a Poisson INAR(1) model and data generated from a different model, namely a Poisson INARCH(1) model,  provided that the series is long enough. The opposite situation, discussed later, is plotted in the right panel.

To sum up, the presented results show that the power of the test is satisfactory provided that the alternative is far enough from the null hypothesis, or the sample size is large.
Moreover, the simulations reveal some differences in powers of the test for different values of the weight parameter ~$a$. However, it seems that it is not possible to generally recommend one particular value of $a$.  For instance, for the case (c) (innovations following a mixture of a Dirac measure at 0 and Poisson distribution), it seems that $a=0$ performs the best. On the other hand, one can observe the opposite or no effect of $a$ in the remaining settings.

\medskip
The power of the test for the null hypothesis of a Poisson INARCH(1) was studied under the following alternatives: (a) a Negative Binomial INARCH(1), (b) a Poisson INAR(1), (c) a Poisson INGARCH(1,1), and (d) a Poisson INARCH(1) with a level change.

Altervative (a) corresponds to the model~\eqref{inarch} with $G$ being the Negative Binomial distribution with a dispersion parameter $r$, see \cite{zhu2011}. Results for dispersion parameter $r=10$ and $r=40$ are plotted in Figure~\ref{fig-inarch-nb}.  Similarly to INAR(1) with Negative Binomial innovations,  the power of the test decreases with increasing $r$, and  in this case, larger values of $a$ seem to lead slightly larger power compared to $a=0$.

Alternative (b) corresponds to a Poisson INAR(1),  and it is presented in  Figure~\ref{fig-inar-inarch}, right panel.  In this case, it seems that larger values of $a$ might yield noticeable larger power compared to $a=0$. Moreover, the effect of $a$ is more substantial here compared to the  situation in the left panel, i.e.~testing of INAR(1) under INARCH(1) alternative.

In alternative (c),  the data are generated from a Poisson  INGARCH(1,1) model of the form
\[ Y_t|\mathcal{I}_t \sim \mathsf{Po}(\lambda_t), \ \lambda_t=\theta_1+\theta_2 Y_{t-1} + \delta \lambda_{t-1},\] see e.g. \cite{ferland}. Figure~\ref{fig-inarch-garch} shows results for $\theta_1=0.1$ and $(\theta_2,\delta)=(0.45,0.5)$ and $(\theta_2,\delta)=(0.25,0.70)$. The unconditional mean of $Y_t$ is equal to $2$ in both settings. 
Clearly, the power growths with an increasing sample size $T$. When comparing different values of $a$, it seems that for larger sample sizes ($T=500$) larger values of $a$ lead to a larger power compared to $a=0$.  However, for smaller sample sizes ($T=100$) one might observe the opposite. Furthermore, when keeping $\theta_1$ and the unconditional mean fixed (i.e.~keeping $\theta_2+\delta$ fixed), the power does not always grow with increasing $\delta$ for a fixed sample size $T$. Specifically, we do observe this behavior for $T=500$, but opposite results are obtained for $T=100$ and $T=250$.

Finally, alternative (d) considers a Poisson INARCH(1) model with a level shift, i.e.~a model of the form
\[Y_t|\mathcal{I}_{t-1}\sim \mathsf{Po}(\lambda_t), \ \lambda_t= \theta_1+\theta_2Y_{t-1} + \delta I[t\geq \tau_0],\]
where $\tau_0$ is a specific time moment, and $I[\cdot]$ is an indicator function.
This model is a special case of a Poisson INARCH(1) with an intervention in \cite{fokianos-fried}.
Various choices for the parameters $\theta_1,\theta_2$, and $\delta$ where considered, together with $\tau_0=\phi \cdot T$, with $\phi\in(0,1)$. Results for $(\theta_1,\theta_2)=(4,0.60)$ and $\delta=4$ and $\delta=8$ are presented in Figure~\ref{fig-inarch-ls}.  For simplicity, we restrict to $a=1$, $T=50$ and $100$ and $\phi=0.1$ and  $\phi=0.3$.
As expected, the power of the test growths substantially with the size of the shift $\delta$. For a fixed value of $\phi$, the power also slowly increases with an increasing sample size $T$.  Furthermore, the simulations  indicate that  different values of $a$ seem to lead to  comparable powers and no general recommendation about ``the most appropriate'' value for $a$ can be given here (results not shown).

\subsection{{\rm{{\bf{INAR(2)}}}}}

In order to illustrate the behavior of the bootstrap test for higher order models, we  present a short simulation study for a Poisson INAR(2) model with $p_1=0.3$, $p_2=0.2$ and $\theta=\mathbb{E} \varepsilon_t=5$.  Under the alternatives, we consider INAR(2) models with innovations $\eps_t$ following (a) the Negative Binomial distribution with dispersion parameter $r$ and (b) a mixture of Poisson and Dirac measure at 0 with weights $1-\phi$ and $\phi$, respectively. For the sake of brevity only results for $a=0$ are shown.

\begin{figure}[htb]
\centering
\includegraphics[width=0.48\textwidth]{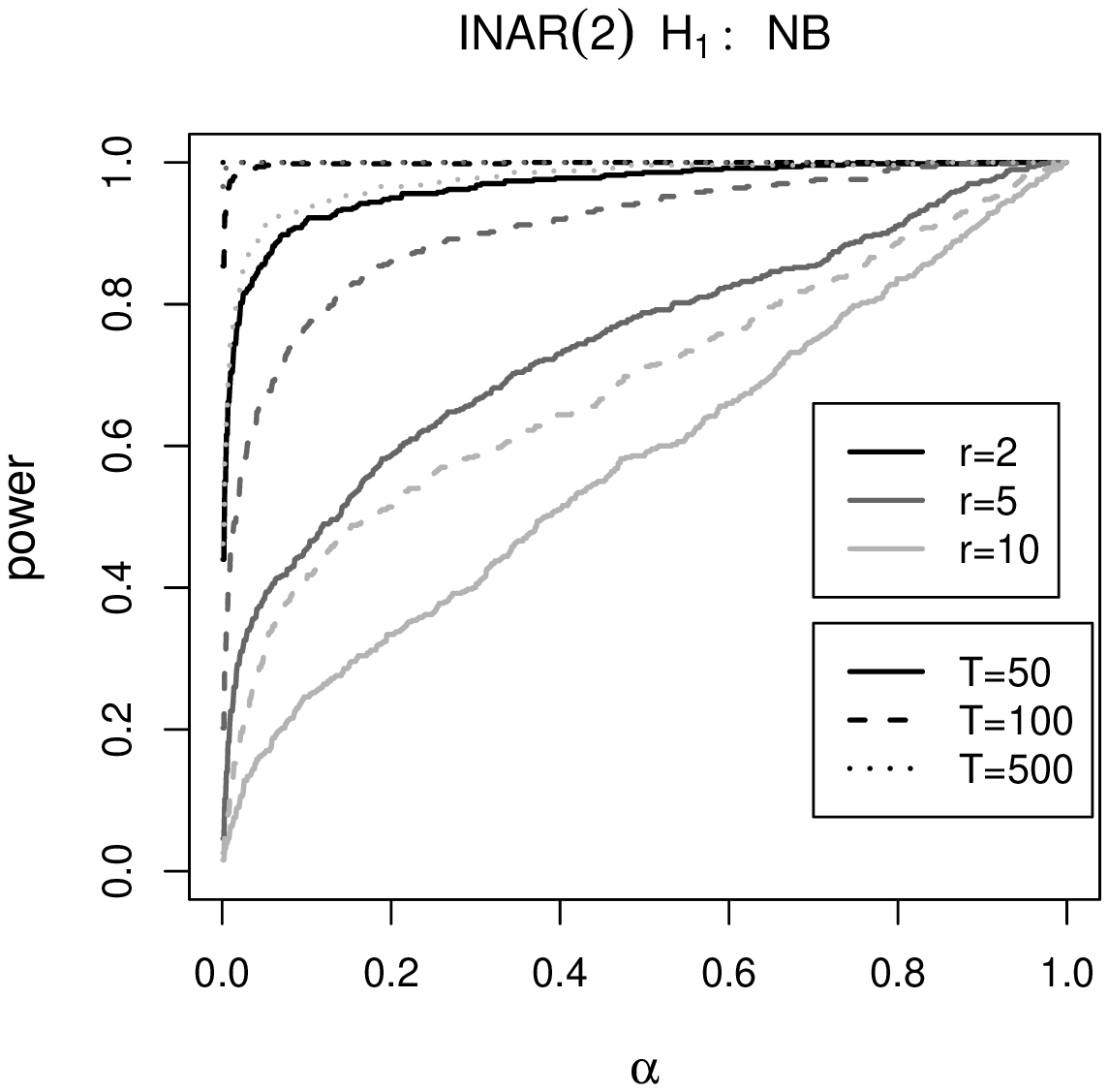}
\includegraphics[width=0.48\textwidth]{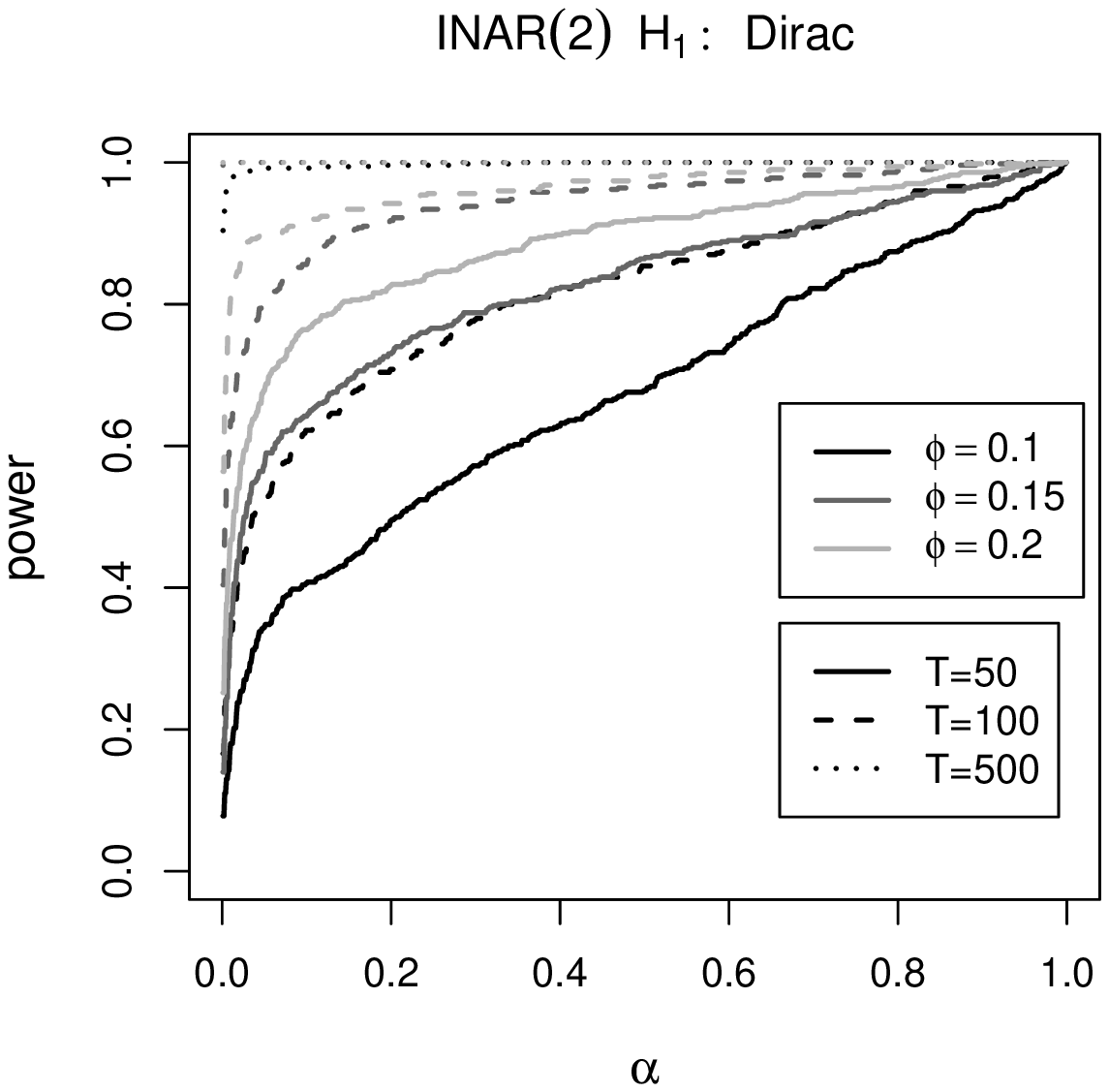}
\caption{Power of the test for a Poisson INAR(2) model under an INAR(2)  with innovations following  (a) the Negative Binomial distribution with dispersion parameter $r$ (left panel) and (b) the mixture of  Poisson and Dirac measure at 0 with weight  $\phi$ (right panel).}\label{f-INAR2dim}
\end{figure}

Table~\ref{t-INAR2dim} indicates that under the null hypothesis, the test generally keeps the prescribed significance level $\alpha$. The power of the test under the alternatives is plotted in Figure~\ref{f-INAR2dim}.
The results are in correspondence with those of INAR(1). In particular, under the alternative (a) the power decreases with increasing dispersion parameter $r$. Under the alternative (b), the power increases with the value of $\phi$. In both cases, the power growths with an increasing sample size $T$. For $T=500$ we get very high power for all the considered settings --- in particular, the power is always greater than 90\% for $\alpha=0.05$.

\begin{table}[htb]
\centering
\tbl{Size of the test for the Poisson INAR(2) model.}{
\begin{tabular}{r|rrr}
\toprule
 &\multicolumn{3}{c}{$\alpha$}\\
$T$  & 0.01 & 0.05 & 0.1 \\
\colrule
50 & 0.006 & 0.038 & 0.096 \\
  100 & 0.010 & 0.040 & 0.076 \\
  250 & 0.010 & 0.044 & 0.106 \\
  500 & 0.006 & 0.050 & 0.094 \\
   \botrule
\end{tabular}}\label{t-INAR2dim}
\end{table}

\section{Application}\label{sec_apl}

We illustrate the goodness-of-fit test on five time series previously analyzed in \cite{freeland}. The data
consist of monthly number of claims of short-term disability benefits made by injured workers to
the British Columbia Workers Compensation Board (WCB). The recorded period is January 1985 to December 1994. The five series correspond to five different injury categories: burn injuries, soft tissue injuries, cuts, dermatitis and dislocations. The first series of burn injuries is corrected by excluding one long duration claimant according to \cite[p.~159]{freeland}.

Freeland \cite{freeland} considered an INAR(1) model for the five series. Except series 3 (cuts) the classical stationary INAR(1) with Poisson marginals was chosen for the analysis. For series 3,
an INAR(1) model with Poisson marginals and seasonality modelled by trigonometric functions was fitted. This series was further investigated e.g. by Zhu and Joe \cite{zhu-joe}, who improved the previous model by considering also Negative Binomial marginals.

We applied the suggested goodness--of--fit test in order to test the appropriateness of the stationary Poisson INAR(1) for the five series. The test statistic $S_T$ was computed using the weight function $w(u)=u^a$ for $a=0,1,2,5,10$. The corresponding p-value was computed from 999 bootstrap samples.

\begin{table}[htb]
\tbl{p--value of the goodness of fit test for the five claim series for $a=0,1,2,5,10$.}
{
\begin{tabular}{rrrrrr}
\toprule
& \multicolumn{5}{c}{$a$}\\
series & 0 & 1 & 2 & 5 & 10 \\
  \colrule
1 & 0.653 & 0.652 & 0.698 & 0.771 & 0.862 \\
  2 & 0.561 & 0.596 & 0.579 & 0.662 & 0.665 \\
  3 & 0.020 & 0.002 & 0.001 & 0.000 & 0.000 \\
  4 & 0.270 & 0.270 & 0.295 & 0.386 & 0.506 \\
  5 & 0.444 & 0.486 & 0.579 & 0.661 & 0.749 \\
 \botrule
\end{tabular}}
\label{tab:real-data}
\end{table}

The obtained results are summarized in Table~\ref{tab:real-data}. We can observe that our results corroborate the results of \cite{freeland}. Specifically, for series 1,2,4, and 5 the null hypothesis that the series follows a Poisson INAR(1) model is not rejected. On the other hand, for series~3 we get significant p-values, which indicate that a simple stationary Poisson INAR(1) model is not appropriate here.

\section*{Acknowledgements}  We thank to both anonymous referees for careful reading and a number of valuable comments which led to
improvement  the presentation of the results. The research of Simos Meintanis was partially supported by grant number 11699 of the Special Account for Research Grants (ELKE) of the National and Kapodistrian University of Athens. The research of Marie Hu\v skov\'a
 was partially supported
 by grant GA\v CR P201/12/1277.  The research of \v S\'arka Hudecov\'a
 was partially supported  by  the Czech Science Foundation project ``DYME – Dynamic Models in Economics'' No.~P402/12/G097.

\appendices

 \section{Proofs}

\noindent {\it Proof of  Theorem 4.1.}
Recall that we assume that
 $\{ Y_t\}_{t\in\mathbb{N}}$ is a sequence of stationary and ergodic variables.
 Hence, the sequence
   $\{ h(Y_t)\}_{t\in\mathbb{N}}$  is stationary and ergodic for any measurable function $h$ such that $\mathbb{E} |h(Y_1)|<\infty$. In particular, this holds for
   $h(y)=u^y-(1+p(u-1)^y)$ for all $u\in[0,1]$.
    
  In the following we drop the index $0$ in $p_0$ and $\theta_0$ whenever it is not
  confusing, and $D$ denotes generic constants.

First of all notice that the test statistic $S_T$ in eqn. (\ref{testINAR})   has  asymptotically the same distribution as
\begin{equation}\label{eq42}
{\cal{S}}_T=\int_0^1 \left(\frac{1}{\sqrt T}\sum_{t=2}^T \widehat Z_t(u)\right)^2 w(u) du,
\end{equation}
where $\widehat Z_t(u)=Z_t(u;\widehat p_T,\widehat \theta_T)$ with $Z_t$ defined in \eqref{Zt}.  Let us
study  the process
 $$
\widehat{\cal{Z}}_T (u)=
   \frac{1}{\sqrt T}\sum_{t=2}^T \widehat Z_t(u) =
   \frac{1}{\sqrt T}\sum_{t=2}^T\Big( u^{Y_t}-(1+\widehat p_T(u-1))^{Y_{t-1}} g_{0,\varepsilon} (u;\widehat \theta_T)\Big)
   $$
    
    for a fixed $u\in[0,1]$. Applying the Taylor expansion we get the decomposition:
    $$
     \widehat{\cal{Z}}_T (u)= J_{0T}(u)+J_{1T}(u)+ J_{2T}(u)+R_{T}(u),
     $$
     where $R_{T}(u)$ is a remainder term  (it does not influence the limit behavior)  and
    \begin{align*}
      J_{0T} (u)=&\ \frac{1}{\sqrt T}\sum_{t=2}^T\Big( u^{Y_t}-(1+p(u-1))^{Y_{t-1}} g_{\varepsilon} (u; \theta)\Big),\\
        J_{1T} (u)= &- \sqrt T(\widehat p_T-p)\frac{1}{ T}\sum_{t=2}^T\Big(  (1+ p(u-1))^{Y_{t-1}-1} Y_{t-1} (u-1)\Big)
    g_{\varepsilon} (u; \theta),\\
     J_{2T} (u)= &-\sqrt T  (\widehat \theta_T-\theta)\frac{1}{ T}\sum_{t=2}^T (1+p(u-1))^{Y_{t-1}}
     \Big( \frac{\partial g_{\varepsilon} (u;\theta)}{\partial
     \theta}\Big).
            \end{align*}
 Standard arguments give that
  \begin{equation*}
|  R_{T}(u)|
    \leq D\Big(T(\widehat p_T-p)^2 + T (\widehat \theta_T-\theta)^2\Big) \Big( \frac{1}{ T^{3/2}} \sum_{t=2}^T
      Y_{t-1}^2 + \frac{1}{\sqrt T}\Big) (1+v(u))
\end{equation*}
      uniformly for $ u\in [0,1]$  for some $D>0$. This together with assumptions (A.3) --  (A.4)   immediately implies
       \[
   \int_0^1R^2_T(u) w(u) du=o_P(1).
   \]
Next we study $J_{1T}(u)$  and $J_{2T}(u)$. In view of the
ergodicity of $\{Y_t\}_{t\in \mathbb N}$ and smoothness
     of $J_{1T}(u)$  and $J_{2T}(u)$ in $u\in[0,1]$ we observe that
       as $T\to \infty$
\begin{align*}
    \frac{1}{T}\sum_{t=2}^T \Big((1+ p(u-1))^{Y_{t-1}-1} Y_{t-1} (u-1)\Big) g_{\varepsilon} (u;\theta)&\to h_1(p, \theta;u), \, a.s.,\\
    \frac{1}{T}\sum_{t=2}^T \Big((1+p(u-1))^{Y_{t-1}}\Big) \frac{\partial g_{\varepsilon} (u;\theta)}{\partial \theta}&\to h_2(p, \theta;u),\, a.s.
  \end{align*}
     for $u\in [0,1]$, where
     \begin{align*}
      h_1(p, \theta;u)=& \frac{\partial  g_Y(1+p(u-1))}{\partial p} (u-1)g_{\varepsilon} (u;\theta)\\
       h_2(p, \theta;u)=&   \mathbb{E}\Big((1+p(u-1))^{Y_{t-1}}\Big) \frac{\partial g_{\varepsilon} (u;\theta)}{\partial \theta}
       = g_Y(1+p(u-1)) \frac{\partial g_{\varepsilon} (u,\theta))}{\partial \theta}.
        \end{align*}
    Since
    $$
    \Big|\Big((1+ p(u-1))^{Y_{t-1}-1} Y_{t-1} (u-1)\Big) g_{\varepsilon} (u;\theta)\Big|\leq  Y_t,\quad  u\in[0,1],
         $$
         and
         $$
          \Big|\Big((1+p(u-1))^{Y_{t-1}}\Big) \frac{\partial g_{\varepsilon} (u;\theta)}{\partial \theta}\Big|\leq  \Big|\frac{\partial g_{\varepsilon} (u;\theta)}{\partial \theta}\Big|,\quad  u\in[0,1]
          $$
          and since the random variables on the r.h.s. have finite  second moments, the uniform ergodicity theorem, see \cite[Theorem 6.2]{rao}, can be applied  and it
      further implies that
               \[
               \int_0^1 \Big(J_{1T}(u)-J_{1T0}(u)\Big)^2 w(u)du +
            \int_0^1 \Big(J_{2T}(u)-J_{2T0}(u)\Big)^2 w(u)du=o_P(1),
                 \]
                  where
                     \begin{align*}
           J_{1T0}(u)  &= - \frac{1}{\sqrt T}\sum_{t=2}^T            \ell_{1} ({\bf{Y}}_{t-q}; p, \theta)
         h_1(p, \theta;u),\\
       J_{2T0}(u)  &= -\frac{1}{\sqrt T}\sum_{t=2}^T \ell_2 ({\bf{Y}}_{t-q}; p,\theta)
       h_2(p, \theta;u).
       \end{align*}
      Hence,    it remains to study the asymptotic behavior of
        $$
         Q_{1T}=
                \int_0^1 \Big( J_{0T}(u)+ J_{1T0}(u)+J_{2T0}(u)\Big)^2 w(u)du.
                $$

In order to obtain the limiting distribution of $\cal{S}_T$ in \eqref{eq42} we apply  Theorem
22 in \cite[p.~380-381]{ibragimov}. We  need to
verify its assumptions. Particularly, we need to show that
  \begin{enumerate}
  \item[(I.)] $\sup_T \mathbb{E} Q_{1T}<\infty,$

  \item[(II.)] $\mathbb{E}|J^2_{0T}(u)-J^2_{0T}(s)|+\mathbb{E}|J^2_{1T0}(u)-J^2_{1T0}(s)|+
  \mathbb{E}|J^2_{2T0}(u)-J^2_{2T0}(s)|\leq D |u-s|^{\kappa},\,$  for some $\kappa>0$ and all
  $u,s\in[0,1],$

  \item[(III.)]  for any  real numbers $a_{vj},\, u_{v}\in [0,1],\, v=1,\ldots,k,\, j=0,1,2 $
   the asymptotic distribution of $\sum_{v=1}^k (a_{0v}J_{0T}(u_{v})+ a_{1v} J_{1T0}(u_{v})+ a_{2v}J_{2T0}(u_{v}))$  is normal with zero mean.
   \end{enumerate}

    To check  validity of (I.)   we notice that  $J_{0T}(u)+ J_{1T0}(u)+J_{2T0}(u)$ are   sums of martingale differences for  any fixed $u$.   Thus,  direct calculations  give
   \begin{align*}
   \mathbb{E}&(J_{0T}(u)+ J_{1T0}(u)+J_{2T0}(u))^2=\frac{T-1}{T} \mathbb{E}\Big(
    u^{Y_t}-(1+p(u-1))^{Y_{t-1}} g_{\varepsilon} (u;\theta)\\
    &-\ell_{1} ({\bf{Y}}_{t-q}; p, \theta)h_1(p, \theta;u)- \ell_{2} ({\bf{Y}}_{t-q}; p, \theta)h_2(p, \theta;u)\Big)^2
   \\
   &\leq D \Big( g^2_{\varepsilon} (u;\theta)+ h^2_1(p, \theta;u)\mathbb{E} \ell^2_{1} ({\bf Y}_{t-q}; p, \theta)
   + h^2_2(p, \theta;u)
   \mathbb{E}\ell^2_{2} ({\bf Y}_{t-q}; p, \theta) +1\Big).
   \end{align*}
   This easily implies (I.).  In addition,  a central limit theorem for sums  of martingale differences can be applied here and this
     immediately gives (III.).

  Concerning  (II.)  notice that
  \begin{align*}
   \mathbb{E}|J^2_{0T}(u_1)-J^2_{0T}(u_2)|&\leq \Big( \mathbb{E}(J_{0T}(u_1)- J_{0T}(u_2))^2
    \mathbb{E}(J_{0T}(u_1)+ J_{0T}(u_2))^2\Big)^{1/2}\\
    &\leq D \Big( \mathbb{E}(J_{0T}(u_1)- J_{0T}(u_2))^2\Big)^{1/2},
       \end{align*}
       where we used the Cauchy-Schwarz inequality and   (I.).
       Since $J_{0T}(u_1)- J_{0T}(u_2)$ is  the sum of a martingale differences sequence  we also get that
       \begin{align*}
    & \mathbb{E}(J_{0T}(u_1)- J_{0T}(u_2))^2= \frac{T-1}{T}\mathbb{E}( u_1^{Y_t} -u_2^{Y_t}
    \\&- \big((1+ p(u_1-1))^{Y_{t-1}}
    - (1+ p(u_2-1))^{Y_{t-1}}\big)g_{\varepsilon} (u_1;\theta)
    \\&-
    (1+ p(u_2-1))^{Y_{t-1}} (g_{\varepsilon} (u_1;\theta)-g_{\varepsilon} (u_2;\theta))\Big)^2
   \\&
    \leq  D(|u_1-u_2|^2 \mathbb{E} Y_t^2 +  |g_{\varepsilon} (u_1;\theta)-g_{\varepsilon} (u_2;\theta)|)\leq D |u_1-u_2|,
    \end{align*}
  where we used the following simple inequalities:
    \begin{align*}
    & | u_1^{Y_t} -u_2^{Y_t}|\leq |u_1-u_2| Y_t(u_1^{Y_t-1}+u_1^{Y_t-1})\leq 2 |u_1-u_2| Y_t,\\
       &|  (1+ p(u_1-1))^{Y_{t-1}}
    - (1+ p(u_2-1))^{Y_{t-1}}|
    \leq D  |u_1-u_2| Y_{t-1}\\
    &  |g_{\varepsilon} (u_1,\theta)-g_{\varepsilon} (u_2,\theta)| \leq D |u_1-u_2|
        \end{align*}
      for some $D>0$. The last inequality follows from the assumption (A.1).

    Concerning $J_{jT0}(u_1)- J_{jT0}(u_2),j=1,2,$  we see that
       $$
 \mathbb{E}|J_{jT0}(u_1)- J_{jT0}(u_2)|^2\leq D \mathbb{E}\ell^2_j({\bf Y}_{t-q};p, \theta)
 |h_j(p,\theta;u_1)-h_j(p,\theta;u_2)|^2, \, j=1,2.
 $$
  By assumptions $\mathbb{E}\ell^2_j({\bf Y}_{t-q};p, \theta)$ is finite and does not depend on $u_1, u_2$. Hence, we just need that
      \[
      |h_j(p,\theta;u_1)-h_j(p,\theta;u_2)|^2\leq D |u_1-u_2|^{\kappa_0},\quad j=1,2,
      \]
       for some $\kappa_0>0$ which holds  true by the assumptions. \qed

   \medskip

   \noindent {\it Proof of  Theorem 4.4.} The proof follows the same lines as
   Theorem~\ref{thm_2_2}. 
    Clearly,
   \begin{align*}
     \exp & \{(\widehat \theta_{1T} +\widehat \theta_{2T} Y_{t-1})(u-1)\}=
     \exp\{( \theta_1 +\theta_2 Y_{t-1})(u-1)\}\\
      & \times \Big( 1+
     ((\widehat \theta_{1T}-\theta_1) +(\widehat \theta_{2T} -\theta_2)Y_{t-1})(u-1) \Big) +R_T(u),\\
     |R_{t,T}(u)|&\leq   \exp\{( \theta_1 +\theta_2 Y_{t-1})(u-1)\}
   \Big(  (\widehat \theta_{1T}-\theta_1)^2 +(\widehat \theta_{2T} -\theta_2)^2Y^2_{t-1}\Big).
 \end{align*}
      By the assumptions
     \begin{align*}
     \int_0^1\frac{1}{T} \large(\sum_{t=2}^T R_{t,T}(u)\large)^2 w(u) du&=O_P\Big(\frac{1}{T}
     \int_0^1  \mathbb{E}\exp\{2
     ( \theta_1 +\theta_2 Y_{t-1})(u-1)\} (1+Y_{t-1}^4) w(u)du\Big)\\
     &=O_P(T^{-1}).
     \end{align*}
After a few steps similar to the steps in INAR(1) we get that under the considered assumptions  and under the null hypothesis, $S_T$ has the same limit distribution as
  \begin{align*}
   \int_0^1 &\frac{1}{T}\Big(\sum_{i=2}^T \Big(u^{Y_t}- \exp\{(\theta_1+\theta_2Y_{t-1})(u-1)\}(u-1)\\+
   &\ell_1({\bf Y}_{t-q}; \theta_1, \theta_2) r_1(u ;\theta_1,
      \theta_2)\\+&
     \ell_2({\bf Y}_{t-q}; \theta_1,\theta_2) r_2(u;\theta_1,
      \theta_2)\Big)^2 w(u) du.
   \end{align*}
   Hence, the assumptions of Theorem 22 in \cite{ibragimov} are fulfilled and our theorem is proved.  \qed

\end{document}